\documentclass[10pt]{amsart}
\PassOptionsToPackage{shortlabels}{enumitem}
\usepackage{amssymb,enumitem,mathtools}
\usepackage{hyperref,aliascnt,tikz-cd,needspace}
\hypersetup{hidelinks}
\usepackage[english]{babel}
\def\today{\number\day\space\ifcase\month\or January\or February\or
 March\or April\or May\or June\or July\or August\or September\or
 October\or November\or December\fi\space\number\year}

\addto\extrasenglish{}
\theoremstyle{definition}
\newtheorem{lma}{Lemma}[section]

\newaliascnt{thmCt}{lma}
\newtheorem{thm}[thmCt]{Theorem}
\aliascntresetthe{thmCt}

\newaliascnt{corCt}{lma}
\newtheorem{cor}[corCt]{Corollary}
\aliascntresetthe{corCt}

\newaliascnt{propCt}{lma}
\newtheorem{prop}[propCt]{Proposition}
\aliascntresetthe{propCt}

\newaliascnt{dfCt}{lma}
\newtheorem{df}[dfCt]{Definition}
\aliascntresetthe{dfCt}

\newaliascnt{remCt}{lma}
\newtheorem{rem}[remCt]{Remark}
\aliascntresetthe{remCt}

\newaliascnt{notnCt}{lma}
\newtheorem{notn}[notnCt]{Notation}
\aliascntresetthe{notnCt}

\newaliascnt{exCt}{lma}
\newtheorem{ex}[exCt]{Example}
\aliascntresetthe{exCt}

\newcounter{theoremintro}

\newtheorem{thmintro}[theoremintro]{Theorem}

\newaliascnt{corintroCt}{theoremintro}
\newtheorem{corintro}[corintroCt]{Corollary}
\aliascntresetthe{corintroCt}

\numberwithin{equation}{section}

\newcommand{\indicator}{\mathdsOne}
\newcommand{\mathdsOne}{\mathord{\mathrm{1}\mkern-4mu\mathrm{l}}}
\newcommand{\Z}{{\mathbb{Z}}}
\newcommand{\R}{{\mathbb{R}}}
\newcommand{\C}{{\mathbb{C}}}
\newcommand{\N}{{\mathbb{N}}}
\newcommand{\T}{{\mathbb{T}}}
\newcommand{\Tr}{{\mathrm{Tr}}}
\newcommand{\id}{{\mathrm{id}}}
\newcommand{\Aut}{{\mathrm{Aut}}}

\newcommand{\cc}{\simeq_{\mathrm{cc}}}
\newcommand{\uu}[1]{\|#1\|_{2,u}}
\newcommand{\HS}[1]{\|#1\|_{\mathrm{HS}}}

\title[]{Equivariant $\mathcal{Z}$-stability for amenable actions}
\date{\today}

\author[Eusebio Gardella]{Eusebio Gardella}
\address{Eusebio Gardella
Department of Mathematical Sciences, Chalmers University of
Technology and University of Gothenburg, Gothenburg SE-412 96, Sweden.}
\email{gardella@chalmers.se}
\urladdr{www.math.chalmers.se/~gardella}

\begin{document}
\raggedbottom
\begin{abstract}
We prove that every amenable action of a countable discrete group on a separable, simple,
unital, nuclear, $\mathcal{Z}$-stable C*-algebra absorbs the trivial action on $\mathcal{Z}$ up
to cocycle conjugacy. The group may be nonamenable, as in Suzuki's recent examples, and no
assumptions on outerness or on the tracial state space are needed. This extends OpenAI's result for
amenable groups to amenable actions. Amenability of the action cannot be omitted in general, as examples of the author and Lupini show. More
generally, we show that every tracially amenable action on a separable, unital, $\mathcal{Z}$-stable
C*-algebra with traces admits invariant central copies of all matrix algebras, without
simplicity or nuclearity. For amenable actions on simple nuclear algebras, equivariant property (SI) allows one to deduce $\mathcal{Z}$-stability from this. 
As applications, we obtain absorption
of $\mathcal{Z}$ for intermediate subalgebras in the crossed product, and also 
pure infiniteness for outer amenable actions of nonamenable groups.  

The basic techniques involving Gaussian variables originate in OpenAI's 
preprint, and are also the main new tool in its 
proposed solution to the Toms--Winter conjecture. We extend these techniques in order to deal with the
weight functions coming from amenability, and one of the main new ingredients is the fact that one can
replace an amenable action by an auxiliary one for which the amenability weights are central.

In order to make the probabilistic tools accessible to the community beyond what is used here, 
we include an appendix where the theory is developed assuming only an elementary background, and where the relevant estimates are proved in detail. 
\end{abstract}
\maketitle

\vspace{-1cm}

\tableofcontents

\vspace{-1cm}

\section{Introduction}

Tensorial absorption has long played a central role in the structure and classification of
operator algebras. For a von Neumann algebra $M$, absorption of the hyperfinite $\mathrm{II}_1$
factor $\mathcal R$ means that $M\cong M\bar\otimes\mathcal R$. McDuff's characterization of
this property in terms of central sequences \cite{McDuff} connected a global tensor product
decomposition with approximately central matrix algebras. This connection was fundamental in
Connes' classification of injective factors \cite{Connes}.

The corresponding theme in C*-algebras developed through absorption of UHF-algebras, the Cuntz
algebra $\mathcal O_\infty$, and the Jiang--Su algebra $\mathcal{Z}$ \cite{JS,TW}. In
particular, $\mathcal{Z}$-stability, the condition $A\cong A\otimes\mathcal{Z}$, became a basic
regularity property in the classification of simple nuclear C*-algebras; see, for example,
\cite{GLN2,CETWW}. The Jiang--Su algebra has the same $K$-theory and tracial state space as the
complex numbers, while its absorption supplies approximately central structure that need not be
visible in these invariants.

For group actions, the analogous question is whether this regularity can be obtained
equivariantly. In von Neumann algebras, absorption of the trivial action on $\mathcal R$
already played an important role in Ocneanu's classification of outer actions of countable
amenable groups on $\mathcal R$ \cite{Ocneanu}; see also \cite{SWMcDuff}. In C*-dynamics, one
asks when an action on a $\mathcal{Z}$-stable algebra absorbs the trivial action on $\mathcal{Z}$
up to cocycle conjugacy.

Equivariant absorption and regularity of crossed products have been studied by several methods. Sato
established results for automorphisms of $\mathcal{Z}$ \cite{SatoRohlin}, while Matui--Sato
treated strongly outer actions \cite{MS1,MS2}. Sato \cite{SatoActions} and
Gardella--Geffen--Naryshkin--Vaccaro \cite{GGNV} also obtained $\mathcal{Z}$-stability of
crossed products without first establishing equivariant absorption. The latter result used
dynamical comparison, while Kerr obtained $\mathcal{Z}$-stability for crossed products of free
minimal almost finite actions \cite{Kerr}.

For amenable acting groups, Szab\'o established equivariant absorption without an outerness
assumption when there were finitely many extreme tracial states \cite{SI1}.
Gardella--Hirshberg--Vaccaro treated finite dimensional Bauer boundaries with restrictions on
trace orbits \cite{GHV}, and Wouters removed the orbit restrictions for single automorphisms
\cite{Wouters}. Szab\'o--Wouters identified local equivariant property Gamma as the obstruction
to equivariant $\mathcal{Z}$-stability for actions of amenable groups \cite{SWGamma}. The recent
preprint \cite{OAI} removed the remaining restrictions on traces for such actions on simple,
separable, unital, nuclear, stably finite, $\mathcal{Z}$-stable C*-algebras. Its new idea was to
use Gaussian quadratic forms to construct invariant central trace splittings, with estimates
independent of the geometry of the tracial state space.

In this paper we show that \emph{every amenable action} on a separable, simple, unital,
nuclear, $\mathcal{Z}$-stable C*-algebra is equivariantly $\mathcal{Z}$-stable. The acting group
need not be amenable, and we impose no conditions on outerness, the tracial state space, or the
action on traces. Thus the answer is the same for the amenable actions of nonamenable groups
constructed by Suzuki \cite{Suzuki} as for actions of amenable groups. We write $\cc$ for
cocycle conjugacy.

\begin{thmintro}\label{mainThmA}
Let $A$ be a separable, simple, unital, nuclear, $\mathcal{Z}$-stable C*-algebra. Let $G$ be a
countable discrete group, and let $\alpha\colon G\to\Aut(A)$ be an amenable action. Then
\[
 (A,\alpha)\cc(A\otimes\mathcal{Z},\alpha\otimes\id_{\mathcal{Z}}).
\]
\end{thmintro}

Suzuki's actions of nonabelian free groups on simple unital AF-algebras \cite[Theorem 1.3 and
Remark 1.5]{Suzuki} give concrete examples to which \autoref{mainThmA} applies. Such actions
have no invariant tracial state, since an amenable action on a unital C*-algebra with an
invariant tracial state forces the group to be amenable \cite[Lemma 2.9]{GGKNV}. The
restriction to amenable actions is meaningful: in \cite{GL}, the author and Lupini showed that the Bernoulli shift
on a finite strongly self-absorbing C*-algebra absorbs the trivial action on $\mathcal{Z}$ if
and only if the acting group is amenable, which is equivalent to the Bernoulli shift itself being amenable.

The extension from amenable groups to amenable actions is far from 
immediate. Scalar F\o lner averages are not available for nonamenable groups, and crucially, the weights
supplied by amenability of the action take values in an algebra and need not be central.
Moreover, the tracial averaging used to combine local constructions cannot be based on an
invariant trace, because such traces need not exist. Both the construction of trace splittings and the
passage from local solutions to invariant matrix algebras therefore require new arguments.

We first construct a commutative $G$-invariant algebra in a tracial relative commutant, and
pass to an auxiliary action for which the amenability weights are central; this is one of the main new ingredients and is likely to be useful in other contexts as well. Under an extreme
trace these weights become scalars, but both their values and the covariances of the Gaussian
variables vary with the trace. We prove joint moment estimates and a central-weight version of
\cite[Proposition 4.1]{OAI}, uniform in both sets of parameters. A new averaging lemma for
central weights then makes the tracial local-to-global principle of \cite{SWGamma} available
for the auxiliary action. The resulting invariant matrix algebras can be transferred back to the
original action.

The tracial result is more general than the setting of \autoref{mainThmA}: it requires neither
simplicity nor nuclearity, and assumes only tracial amenability in the sense of \cite{GGKNV}.
Fix a free ultrafilter $\omega$ on $\N$. We write $T(A)$ for the tracial state space of $A$ and
$A^\omega$ for the uniform tracial ultrapower, $\alpha^\omega$ for the induced action, and
$A^\omega\cap A'$ for the relative commutant of the constant image of $A$. These notions are
recalled in \autoref{sec:prelim}.

\begin{thmintro}\label{mainThmB}
Let $A$ be a separable, unital, $\mathcal{Z}$-stable C*-algebra with $T(A)\ne\emptyset$. Let
$G$ be a countable discrete group, and let $\alpha\colon G\to\Aut(A)$ be tracially amenable.
For every integer $m\ge2$, there is a unital homomorphism
\[
 M_m\to (A^\omega\cap A')^{\alpha^\omega}.
\]
\end{thmintro}

This is an equivariant McDuff-type conclusion for arbitrary countable groups and without
restrictions on trace orbits. When $T(A)\ne\emptyset$, equivariant property (SI) gives
\autoref{mainThmA} from \autoref{mainThmB}; it holds under the hypotheses of \autoref{mainThmA}
by \cite[Corollary 3.3]{SI}. This last step explains why amenability of the action is required
in \autoref{mainThmA}: tracially amenable actions on classifiable algebras need not have
equivariant property (SI) \cite[Example A]{GKV}. The case $T(A)=\emptyset$ follows from known
results, as explained in the proof of \autoref{thm:absorption}.

The equivariant conclusion also gives, for example, regularity of crossed products and all
their intermediate algebras. We write $A\rtimes_{\alpha,r}G$ for the reduced crossed product.

\begin{corintro}\label{cor:intro-intermediate}
Let $A$ be a separable, simple, unital, nuclear, $\mathcal{Z}$-stable C*-algebra, let $G$ be a
countable discrete group, and let $\alpha\colon G\to\Aut(A)$ be amenable. Then every
intermediate C*-algebra
\[
 A\subseteq C\subseteq A\rtimes_{\alpha,r}G
\]
is $\mathcal{Z}$-stable. In particular, $A\rtimes_{\alpha,r}G$ is $\mathcal{Z}$-stable.
\end{corintro}

For nonamenable groups, this gives pure infiniteness when the action is outer. Recall that
$\alpha$ is \emph{outer} if $\alpha_g$ is outer for every $g\ne e$, and that a \emph{Kirchberg
algebra} is a separable, simple, nuclear, purely infinite C*-algebra.

\begin{corintro}\label{cor:intro-pure}
Let $A$ be a separable, simple, unital, nuclear, $\mathcal{Z}$-stable C*-algebra, let $G$ be a
countable discrete nonamenable group, and let $\alpha\colon G\to\Aut(A)$ be amenable and outer.
Then $A\rtimes_{\alpha,r}G$ is a unital Kirchberg algebra.
\end{corintro}

This answers the question motivating \cite{GGKNV} for arbitrary countable nonamenable groups
under the hypotheses above. For groups with weak paradoxical towers, including those containing
a nonabelian free subgroup, pure infiniteness was obtained there under the weaker assumption of
tracial amenability, without nuclearity \cite[Corollary 3.15]{GGKNV}. Our intermediate-algebra
conclusion follows from absorption results for inclusions in \cite{Sarkowicz}; we also give the
central sequence argument in \autoref{sec:applications}.

A further purpose of the paper is to make the Gaussian methods introduced in \cite{OAI}
accessible to operator algebraists. These methods also belong to a broader development:
Gaussian central-limit constructions feature in OpenAI's proposed proof of the unital
Toms--Winter conjecture \cite{OAITW}. We do not use that preprint or assume the correctness of
its claims. The probabilistic arguments in \cite{OAI} are difficult to extract from the
surrounding proof, and the connection with C*-dynamics is initially surprising. The appendix
develops the needed background starting from elementary probability, proves the moment and
Gaussian quadratic-form estimates, and gives a proof of the version of \cite[Proposition
4.1]{OAI} used here. Thus no assertion from \cite{OAI} enters our proof as an unverified black
box. Before beginning the main argument, we explain the connection informally.

\subsection{Gaussian variables and equivariant absorption}\label{sec:intro-gaussian}

A trace turns a selfadjoint element into something that behaves like a real random variable:
its spectral values have a distribution, and the trace of its powers gives the moments of that
distribution. This can be seen without probability theory. Set
$d=\operatorname{diag}(1,-1)\in M_2$, let $\tau_N$ be the normalized trace on
$M_2^{\otimes N}$, and denote by $x_r$ the copy of $d$ in the $r$th tensor factor. Set
$z_N=\frac{1}{\sqrt N}\sum_{r=1}^N x_r$. The tensor product basis is indexed by sign vectors in
$\{-1,1\}^N$, and $z_N$ has diagonal entry $\frac{1}{\sqrt N}\sum_{r=1}^N\varepsilon_r$ at
$(\varepsilon_1,\ldots,\varepsilon_N)$. Since $\tau_N$ assigns weight $2^{-N}$ to every basis
vector, the spectral distribution of $z_N$ under $\tau_N$ is precisely that of a normalized sum
of independent signs. In particular,
\begin{equation}\label{eq:intro-matrix-moments}
 \tau_N(z_N^2)=1
 \quad\text{and}\quad
 \tau_N(z_N^4)=3-\frac{2}{N}.
\end{equation}
The limits are the second and fourth moments of a standard real Gaussian variable. In fact, the
spectral distributions of the $z_N$ converge to the Gaussian distribution.

This observation is useful because the desired invariant central projections should divide all
traces equally. If one could produce a central selfadjoint element fixed by the action, with
the same atomless spectral distribution under every trace, then cutting its spectrum into sets
of equal measure would do precisely that. In a tracial ultrapower the cuts are obtained as
limits of continuous functions. The problem is therefore to make the element nearly invariant
while keeping its spectral distribution spread out, uniformly over the traces.

Absorption of $\mathcal{Z}$ supplies many central copies from which to form normalized sums,
much as the tensor factors did in the example. The group translates of an element in a given
copy may be correlated, and need not commute; different copies, however, can be arranged to
commute. Under an extreme trace, the copies are independent in the sense that the trace of a
product from distinct copies is the product of their traces. Adding centered elements from many
such copies with a factor $1/\sqrt N$ leads to Gaussian joint moments. Here ``centered'' just
means that each element has trace zero. The covariance of the resulting Gaussian variables
records how the action relates the translates, rather than assuming that the translates are
independent.

There remains a difficulty which an ordinary average cannot solve. Even with scalar F\o lner
weights, averaging $N$ independent centered variables with a factor $1/N$ makes their variance
tend to zero, and an almost scalar element has no useful spectrum to divide. Gaussian quadratic
expressions allow the averaging and the preservation of fluctuations to coexist. The weights
and the correlations determine a positive finite dimensional operator $P$ of trace one; if its
eigenvalues are $\lambda_1,\ldots,\lambda_d$ and $\eta_1,\ldots,\eta_d$ are independent
standard Gaussians, the corresponding centered expression is
\[
 Q=\sum_{j=1}^d\sqrt{\lambda_j}(\eta_j^2-1).
\]
Its variance is always $2$, regardless of the correlations and the number of directions. More
is true: the Fourier estimates show that a suitable exponential of $Q$ spreads its spectral
distribution around the circle. In the ultrapower this gives a unitary with Haar spectral
distribution, whose spectrum can be divided into equal arcs. Since the bounds depend only on
the trace-one normalization of $P$, they remain valid as the trace and its observed
correlations vary.

For amenable actions of nonamenable groups, we replace the scalar F\o lner weights by central
weights in an auxiliary algebra. An extreme trace evaluates these as scalars, so the same
picture applies even though the weights vary with the trace. The rest of the paper makes this
argument precise, and the appendix can be read independently as an introduction to the Gaussian
estimates needed for this method.

\subsection*{AI statement}
This work was carried out with the interactive assistance of ChatGPT. It helped
the author understand parts of the preprint \cite{OAI} and suggested proofs of some results. The author checked the arguments used here, substantially rewrote many of them, and organized the material in the way it is presented now.  
The author takes full responsibility for
the mathematical claims and for any inaccuracies or gaps that remain.

\section{Preliminaries}\label{sec:prelim}

In this section, we collect the definitions and results that will be used in the proofs of
\autoref{mainThmA} and \autoref{mainThmB}.
Most of the argument takes place in uniform tracial ultrapowers, where we will construct
invariant matrix algebras commuting with a prescribed separable subset.
We begin by recalling these ultrapowers and the approximation arguments that will allow us to
work in them, and then discuss amenability and local equivariant property Gamma.
The latter is the property needed to pass from constructions in individual tracial
representations to a construction that is uniform over all traces.
We then recall the two absorption results used in the proof: McDuff absorption provides the
matrix algebras in the tracial representations, while equivariant property (SI) allows us to
deduce equivariant $\mathcal{Z}$-stability from \autoref{mainThmB}.
Finally, \autoref{sec:probability} recalls the probability notation used to obtain local
equivariant property Gamma, with references to the development and proofs in
\autoref{sec:probability-appendix}.

All tensor products of C*-algebras are minimal, and all homomorphisms between C*-algebras are
understood to preserve adjoints.
We write $\mathcal{Z}$ for the Jiang--Su algebra, $\tau_{\mathcal{Z}}$ for its unique tracial
state, and $\mathcal R$ for the hyperfinite $\mathrm{II}_1$ factor.
For a unital C*-algebra $B$, we denote by $B_1$ its closed unit ball and by $T(B)$ its tracial
state space, and we keep the free ultrafilter $\omega$ on $\N$ fixed throughout.
A subset of $\N$ is called \emph{$\omega$-large} if it belongs to $\omega$.

\subsection{Uniform tracial ultrapowers}\label{sec:tracial-ultrapowers}
The uniform tracial ultrapower allows us to express approximation in all tracial $2$-seminorms
simultaneously as exact relations, and it is crucial in our work.
We recall its construction next.

\begin{df}\label{def:tracial-ultrapower}
Let $B$ be a unital C*-algebra with $T(B)\ne\emptyset$, and for $b\in B$ and $\tau\in T(B)$ set
\[
 \|b\|_{2,\tau}=\tau(b^*b)^{\frac{1}{2}}
 \quad\text{and}\quad
 \uu b=\sup_{\tau\in T(B)}\|b\|_{2,\tau}.
\]
We denote by $I_{2,\omega}$ the ideal of $\ell^\infty(\N,B)$ given by
\[
 I_{2,\omega}=\bigl\{(b_n)_{n\in\N}\in\ell^\infty(\N,B)\colon\lim_{n\to\omega}\uu{b_n}=0\bigr\},
\]
and call the quotient $B^\omega=\ell^\infty(\N,B)/I_{2,\omega}$ the \emph{uniform tracial ultrapower} of $B$.
\end{df}

\begin{rem}
We denote by $\Delta_B\colon B\to B^\omega$ the homomorphism given by constant sequences.
This homomorphism need not be injective, so expressions such as $B^\omega\cap B'$ always refer
to the relative commutant of its image.
\end{rem}

\begin{notn}
If $\beta\colon G\to\Aut(B)$ is an action, then each $\beta_g$ preserves the uniform
$2$-seminorm, and hence there is an induced action $\beta^\omega\colon G\to\Aut(B^\omega)$
satisfying
\[
 \beta_g^\omega([(b_n)_{n\in\N}])=[(\beta_g(b_n))_{n\in\N}]
\]\nopagebreak
for all $g\in G$ and all bounded sequences $(b_n)_{n\in\N}$ in $B$.
\end{notn}

The traces that we will use on $B^\omega$ are obtained from sequences of traces on $B$, and we define these next.
\begin{df}\label{def:limit-traces}
Given a sequence $(\tau_n)_{n\in\N}$ in $T(B)$, define $\tau\colon B^\omega\to\C$ by
\[
 \tau([(b_n)_{n\in\N}])=\lim_{n\to\omega}\tau_n(b_n)
\]\nopagebreak
for all $(b_n)_{n\in\N}\in\ell^\infty(\N,B)$.
This is a tracial state, called a \emph{limit trace}, and we denote the set of all such traces by $T_\omega(B)$.
For $x\in B^\omega$, we again write $\uu x=\sup_{\tau\in T_\omega(B)}\|x\|_{2,\tau}$.
\end{df}

\begin{rem}\label{rem:ultranorm}
For every bounded sequence $(x_n)_{n\in\N}$ in $B$, we have
\begin{equation}\label{eq:ultranorm}
 \uu{[(x_n)_{n\in\N}]}=\lim_{n\to\omega}\uu{x_n}.
\end{equation}
Indeed, each sequence of traces gives the inequality from left to right.
For the reverse inequality, choose $\tau_n\in T(B)$ such that
$\|x_n\|_{2,\tau_n}>\uu{x_n}-\frac{1}{n}$ for every $n\in\N$.
The limit trace defined by $(\tau_n)_{n\in\N}$ attains the right hand side of
\eqref{eq:ultranorm}, since square roots commute with ultralimits of bounded nonnegative
sequences.
In particular, $T_\omega(B)$ separates the positive elements of $B^\omega$.
The displayed formula allows us to pass between an element of the ultrapower and a bounded representative.
We will also use the elementary inequalities
\begin{equation}\label{eq:2bounds}
 \uu{xy}\le\|x\|\uu y,
 \quad\uu{xy}\le\uu x\|y\|,
 \quad\text{and}\quad
 |\tau(xy)|\le\|x\|_{2,\tau}\|y\|_{2,\tau}
\end{equation}\nopagebreak
for all $x,y\in B^\omega$ and all $\tau\in T_\omega(B)$.
In particular, multiplication is jointly continuous in the uniform $2$-norm on norm bounded
sets, and commutation with a fixed element is preserved under bounded limits in this norm.
\end{rem}

Our constructions will first satisfy finitely many of the required relations up to a prescribed error.
To obtain elements satisfying all the relations exactly, we use Kirchberg's $\varepsilon$-test,
which we recall in the following form; see \cite[Lemma 3.1]{KR}.

\begin{lma}\label{lem:diagonal}
Let $(X_n)_{n\in\N}$ be nonempty sets.
For each $j,n\in\N$, let $f_n^{(j)}\colon X_n\to[0,\infty)$ be a function.
Suppose that for every $k\in\N$ there are elements $x_n^{(k)}\in X_n$, for $n\in\N$, such that
\[
 \lim_{n\to\omega}f_n^{(j)}(x_n^{(k)})<\frac{1}{k}
\]\nopagebreak
for all $j=1,\ldots,k$.
Then there are elements $x_n\in X_n$, for $n\in\N$, such that
$\lim_{n\to\omega}f_n^{(j)}(x_n)=0$ for every $j\in\N$.
\end{lma}

In our applications, the sets $X_n$ will consist of contractions, unitaries, or finite tuples
of contractions, and the functions will measure the failure of the relations we want to impose.
Besides uniform $2$-norms of polynomial expressions, these functions will include suprema of trace expressions.
Let $(a_n)_{n\in\N}$ be a bounded sequence in $B$.
The same argument as in \autoref{rem:ultranorm}, now choosing $\nu_n\in T(B)$ with
$|\nu_n(a_n)|>\sup_{\nu\in T(B)}|\nu(a_n)|-\frac{1}{n}$ for every $n\in\N$, gives
\begin{equation}\label{eq:trace-sup}
 \sup_{\tau\in T_\omega(B)}|\tau([(a_n)_{n\in\N}])|
 =\lim_{n\to\omega}\sup_{\nu\in T(B)}|\nu(a_n)|.
\end{equation}
For example, after choosing a norm dense sequence $(b_j)_{j\in\N}$ in $B_1$, we measure
centrality by the functions $f_n^{(j)}(x)=\uu{[x,b_j]}$ for $x\in B_1$.
For a fixed group element $g$, we likewise measure invariance by $\uu{\beta_g(x)-x}$.

\begin{rem}
We will also use the standard consequence of the $\varepsilon$-test that every norm bounded
ball of $B^\omega$ is complete in the uniform $2$-norm.
Together with \autoref{rem:ultranorm}, completeness allows us to form bounded uniform $2$-norm
limits while preserving products, commutation relations, and invariance under $\beta^\omega$.
In \autoref{rem:spectral-projections}, we explain in detail how to use this observation to
construct projections from unitaries with a prescribed spectral distribution.
\end{rem}

\subsection{Amenability and central weights}\label{sec:amenability}
We next recall the approximation properties of amenable and tracially amenable actions.
The distinction between them is relevant to our two main results: tracial amenability suffices
for the constructions in the uniform tracial ultrapower, whereas amenability, via equivariant
property (SI), is needed to obtain results in the norm ultrapower.

\begin{notn}\label{not:hilbert-module}
Let $G$ be a discrete group, let $A$ be a unital C*-algebra, and let $\alpha\colon G\to\Aut(A)$ be an action.
We write $C_c(G,A)$ for the space of finitely supported functions from $G$ to $A$, and equip it
with the $A$-valued inner product
\[
 \langle\xi,\zeta\rangle=\sum_{g\in G}\xi(g)^*\zeta(g)
\]\nopagebreak
for all $\xi,\zeta\in C_c(G,A)$.
Its completion with respect to the norm $\|\xi\|=\|\langle\xi,\xi\rangle\|^{\tfrac{1}{2}}$ is
the Hilbert $A$-module $\ell^2(G,A)$, on which $A$ acts on the left and on the right by
pointwise multiplication.
The action $\alpha$ induces an action $\widetilde\alpha$ on this module, given by
$\widetilde\alpha_h(\xi)(g)=\alpha_h(\xi(h^{-1}g))$ for all $g,h\in G$ and all
$\xi\in\ell^2(G,A)$.
\end{notn}
We will use the characterization given in \cite[Theorem 3.2]{OS} as the definition of amenability.

\begin{df}\label{def:amenability}\label{def:tracial-amenability}
Let $\alpha\colon G\to\Aut(A)$ be an action of a discrete group on a unital C*-algebra.
If $T(A)\ne\emptyset$, for $\zeta\in\ell^2(G,A)$ set
\[
 \|\zeta\|_{2,u}=\sup_{\tau\in T(A)}\tau(\langle\zeta,\zeta\rangle)^{\frac{1}{2}}.
\]
We say that $\alpha$ is \emph{tracially amenable} if $T(A)\ne\emptyset$ and, for every finite
$F\subseteq A$, every finite $K\subseteq G$, and every $\varepsilon>0$, there is
$\xi\in C_c(G,A)$ with $\|\xi\|\le1$ such that
\[
 \|a\xi-\xi a\|_{2,u}<\varepsilon,
 \quad\uu{\langle\xi,\xi\rangle-1}<\varepsilon,
 \quad\text{and}\quad
 \|\widetilde\alpha_h(\xi)-\xi\|_{2,u}<\varepsilon
\]\nopagebreak
for all $a\in F$ and all $h\in K$; see \cite[Definition 2.2]{GGKNV}.
We say that $\alpha$ is \emph{amenable} if the same approximation holds in norm, that is, if
one can choose $\xi$ with $\|\xi\|\le1$ such that
\[
 \|a\xi-\xi a\|<\varepsilon,
 \quad\|\langle\xi,\xi\rangle-1\|<\varepsilon,
 \quad\text{and}\quad
 \|\widetilde\alpha_h(\xi)-\xi\|<\varepsilon
\]\nopagebreak
for all $a\in F$ and all $h\in K$.
\end{df}
Amenability is defined without any assumption on the tracial state space; in particular, it
also applies when $T(A)=\emptyset$.

\begin{rem}
Since $\|\cdot\|_{2,u}\le\|\cdot\|$, every amenable action on an algebra with nonempty tracial
state space is tracially amenable.
For actions on unital, simple, separable, nuclear, stably finite C*-algebras with strict
comparison, the difference between the two notions is precisely equivariant property (SI), as
we recall in \autoref{thm:amenability-si} below.
\end{rem}

For actions on commutative algebras, amenability has a formulation in terms of positive
functions that is particularly convenient for averaging.
Amenability of an action on $C(Y)$ is equivalent to topological amenability of the induced
action on $Y$, and we use the following formulation from \cite[Definition 2.2]{OS}.

\begin{thm}\label{thm:amenability-weights}
Let $Y$ be a compact Hausdorff space, let $G$ be a discrete group, and let
$\gamma\colon G\to\Aut(C(Y))$ be an action.
The action $\gamma$ is amenable if and only if, for every finite $K\subseteq G$ and every
$\varepsilon>0$, there is a finitely supported family $(d_g)_{g\in G}$ in $C(Y)_+$ such that
\begin{equation}\label{eq:weights}
 \sum_{g\in G}d_g=1
 \quad\text{and}\quad
 \Bigl\|\sum_{g\in G}|\gamma_h(d_{h^{-1}g})-d_g|\Bigr\|<\varepsilon
\end{equation}\nopagebreak
for all $h\in K$.
\end{thm}

In the context of \autoref{thm:amenability-weights}, we refer to the elements of
$(d_g)_{g\in G}$ as \emph{amenability weights}.

We will use \autoref{thm:amenability-weights} when $C(Y)$ is a unital $G$-invariant subalgebra
of the center of another C*-algebra.
More precisely, if $\beta\colon G\to\Aut(B)$ is an action and $E\subseteq Z(B)$ is a unital
$G$-invariant subalgebra on which the action is amenable, then we may choose the weights in
$E$.
Their images in $B^\omega$ are still central, since they commute with every bounded sequence in $B$.
This is the situation to which we will reduce \autoref{mainThmB} in \autoref{prop:auxiliary},
and it allows us to use the weights both to construct approximately invariant unitaries and in
the averaging argument of \autoref{sec:averaging}.

\subsection{Local equivariant property Gamma}\label{sec:prelim-gamma}
The first step in the proof of \autoref{mainThmB} is to construct invariant central projections
that divide every trace equally.
For the tracial local-to-global principle used later, we need to do this relative to a
prescribed subset of the ultrapower, and with the same trace splitting after multiplying by any
element of that subset.
We recall the corresponding form of local equivariant property Gamma from \cite{SWGamma}.

\begin{df}\label{def:gamma}
Let $B$ be a unital C*-algebra with $T(B)\ne\emptyset$, and let $\beta\colon G\to\Aut(B)$ be an
action of a countable discrete group.
We say that $\beta$ has \emph{local equivariant property Gamma} if for every uniform $2$-norm
separable subset $S\subseteq B^\omega$ and every integer $m\ge2$, there are pairwise orthogonal
projections $p_1,\ldots,p_m\in(B^\omega)^{\beta^\omega}\cap S'$ such that
\[
 \sum_{j=1}^m p_j=1\quad\text{and}\quad\tau(sp_j)=\frac{1}{m}\tau(s)
\]\nopagebreak
for all $s\in S$, all $\tau\in T_\omega(B)$, and all $j=1,\ldots,m$.
\end{df}

In the definition above, taking $1\in S$ gives $\tau(p_j)=\frac{1}{m}$ for all
$\tau\in T_\omega(B)$ and all $j=1,\ldots,m$.
The condition on $\tau(sp_j)$ is stronger than this equality, and is needed to combine
constructions made in different tracial representations.
Note also that $S$ is allowed to lie in $B^\omega$, rather than only in the constant image of $B$.

\subsection{Tracial representations and McDuff absorption}\label{sec:prelim-mcduff}
To construct matrix subalgebras, we will first work in tracial representations and apply the
absorption theorem for amenable actions on McDuff von Neumann algebras.
An individual trace need not be invariant, however, so the original action need not extend to
the weak closure in its GNS representation.
We avoid this difficulty by considering all translates of the trace at once, using the following representation.

\begin{df}\label{def:orbit}
Let $C$ be a unital C*-algebra, let $G$ be a countable discrete group, let
$\gamma\colon G\to\Aut(C)$ be an action, and let $\tau\in T(C)$.
We write $(\pi_\tau,\mathcal H_\tau)$ for the GNS representation of $\tau$.
Define the \emph{orbit representation}
$\pi_\tau^\gamma\colon C\to\mathcal B(\ell^2(G)\otimes\mathcal H_\tau)$ and the \emph{orbit von
Neumann algebra} $N_\tau^\gamma$ by
\begin{equation}\label{eq:orbit}
 \pi_\tau^\gamma=\bigoplus_{g\in G}\pi_\tau\circ\gamma_{g^{-1}}
 \quad\text{and}\quad
 N_\tau^\gamma=\pi_\tau^\gamma(C)''\subseteq\mathcal B(\ell^2(G)\otimes\mathcal H_\tau).
\end{equation}
\end{df}

\begin{rem}\label{rem:orbit}
For $h\in G$, denote by $U_h$ the unitary on $\ell^2(G)\otimes\mathcal H_\tau$ satisfying
$(U_h\xi)(g)=\xi(h^{-1}g)$ for all $\xi\in\ell^2(G,\mathcal H_\tau)$ and all $g\in G$.
Then
\[
 U_h\pi_\tau^\gamma(c)U_h^*=\pi_\tau^\gamma(\gamma_h(c))
\]\nopagebreak
for all $c\in C$ and all $h\in G$.
It follows that conjugation by $U_h$ restricts to a normal automorphism of $N_\tau^\gamma$ extending $\gamma_h$.

\end{rem}

\begin{rem}\label{rem:orbittrace}
Adopt the assumptions and notation from \autoref{def:orbit}.
To form the tracial ultrapower of $N_\tau^\gamma$, we choose a faithful normal tracial state as follows.
Denote by $e$ the identity of $G$, and choose strictly positive numbers $c_g$, for $g\in G$,
with $\sum_{g\in G}c_g=1$.
Denote by $\overline\tau$ the faithful normal tracial state on $\pi_\tau(C)''$ extending $\tau$.
For $g\in G$, compression to the $g$th summand gives a normal homomorphism
$p_g\colon N_\tau^\gamma\to\pi_\tau(C)''$ satisfying
$p_g(\pi_\tau^\gamma(c))=\pi_\tau(\gamma_{g^{-1}}(c))$ for all $c\in C$.
Set $\widetilde\tau_g=\overline\tau\circ p_g$ and $\widehat\tau=\sum_{g\in G}c_g\widetilde\tau_g$.
In particular, $\widetilde\tau_g(\pi_\tau^\gamma(c))=\tau(\gamma_{g^{-1}}(c))$ for all $c\in C$.
The traces $\widetilde\tau_g$ jointly separate positive elements, so $\widehat\tau$ is a
faithful normal tracial state satisfying
\begin{equation}\label{eq:orbittrace}
 \widehat\tau(\pi_\tau^\gamma(c))=\sum_{g\in G}c_g\tau(\gamma_{g^{-1}}(c))
 \quad\text{and}\quad
 \tau(c^*c)\le \frac{1}{c_e}\widehat\tau(\pi_\tau^\gamma(c^*c))
\end{equation}\nopagebreak
for all $c\in C$.
\end{rem}

\begin{notn}\label{not:vn-ultrapower}
For a finite von Neumann algebra $M$ with a faithful normal tracial state $\tau_0$, we denote
by $M^\omega$ the quotient of $\ell^\infty(\N,M)$ by the ideal of sequences $(x_n)_{n\in\N}$
satisfying $\lim_{n\to\omega}\|x_n\|_{2,\tau_0}=0$.
Recall that on norm bounded subsets of a finite von Neumann algebra, the strong-* topology is
induced by the $2$-norm associated with any faithful normal tracial state.
Therefore, the ideal of bounded sequences that converge to zero along $\omega$ in this $2$-norm
is independent of the chosen faithful normal trace and is preserved by every normal
automorphism.
The action on $N_\tau^\gamma$ therefore induces an action on its tracial ultrapower even when
$\widehat\tau$ is not invariant.
\end{notn}

Recall that a von Neumann algebra $M$ is \emph{McDuff} if $M\cong M\bar\otimes\mathcal R$.
The absorption theorem for amenable actions on McDuff von Neumann algebras gives the invariant
matrix algebras that we will need in the tracial representations.
We record both conclusions in the following statement.

\begin{thm}[{\cite[Theorems 5.4 and 5.5]{SWMcDuff}}]\label{thm:vn-absorption}
Let $M$ be a finite McDuff von Neumann algebra with separable predual, and let
$\gamma\colon G\to\Aut(M)$ be an amenable action of a countable discrete group.
Then $\gamma\cc\gamma\otimes\id_{\mathcal R}$, and there is a unital homomorphism
\begin{equation}\label{eq:vncriterion}
 M_m\to(M^\omega\cap M')^{\gamma^\omega}
\end{equation}\nopagebreak
for every integer $m\ge2$.
\end{thm}

We use the following standard consequence of $\mathcal{Z}$-stability to check McDuffness in our setting.

\begin{lma}\label{lem:gns}
Let $B$ be a separable, unital, $\mathcal{Z}$-stable C*-algebra with $T(B)\ne\emptyset$, and let $\tau\in T(B)$.
Then $\pi_\tau(B)''$ has separable predual and is McDuff.
\end{lma}

\subsection{Equivariant property (SI) and absorption}\label{sec:prelim-si}
The matrix algebras in \autoref{mainThmB} belong to a tracial ultrapower, while equivariant
$\mathcal{Z}$-stability is characterized by homomorphisms into the norm central sequence
algebra.
We next recall how equivariant property (SI) allows us to pass from one to the other.

\begin{notn}\label{not:norm-ultrapower}
Let $A$ be a unital C*-algebra, and let $\alpha\colon G\to\Aut(A)$ be an action of a countable discrete group.
We write
\[
 A_\omega=\ell^\infty(\N,A)/\bigl\{(a_n)_{n\in\N}\colon\lim_{n\to\omega}\|a_n\|=0\bigr\}
\]\nopagebreak
for the norm ultrapower, denote the induced action by $\alpha_\omega$, and use
$A_\omega\cap A'$ for the central sequence algebra.
If $T(A)\ne\emptyset$, denote by $\pi_A\colon A_\omega\to A^\omega$ the quotient map and set $J_A=\ker(\pi_A)$.
We regard limit traces on $A^\omega$ also as traces on $A_\omega$ by composition with $\pi_A$.
\end{notn}

\begin{df}\label{def:si}
Let $A$ be a separable, simple, unital C*-algebra with $T(A)\ne\emptyset$, and let
$\alpha\colon G\to\Aut(A)$ be an action of a countable discrete group.
We say that $\alpha$ has \emph{equivariant property (SI)} if, whenever
$e,f\in (A_\omega\cap A')^{\alpha_\omega}$ are positive contractions satisfying $e\in J_A$ and
\[
 \inf_{k\in\N}\inf_{\tau\in T_\omega(A)}\tau(f^k)>0,
\]
there is a contraction $s\in (A_\omega\cap A')^{\alpha_\omega}$ such that $s^*s=e$ and $fs=s$.
\end{df}
The next result identifies the additional condition needed to recover amenability from tracial amenability.
In particular, it supplies equivariant property (SI) under the hypotheses of \autoref{mainThmA}
when the tracial state space is nonempty.

\begin{thm}[{\cite[Theorem B]{GKV}}]\label{thm:amenability-si}
Let $A$ be a unital, simple, separable, nuclear, stably finite C*-algebra with strict
comparison, and let $\alpha\colon G\to\Aut(A)$ be an action of a countable discrete group.
Then $\alpha$ is amenable if and only if it is tracially amenable and has equivariant property (SI).
\end{thm}

The implication from amenability to equivariant property (SI) also follows from \cite[Corollary
3.3]{SI}, which applies more generally to amenable actions on nonelementary separable simple
nuclear C*-algebras with strict comparison.
The next result explains how we will use this property to pass from the conclusion of
\autoref{mainThmB} to equivariant $\mathcal{Z}$-stability.

\begin{thm}[{\cite[Theorem 5.4]{SWGamma}}]\label{thm:si-absorption}
Let $A$ be a separable, simple, unital C*-algebra with $T(A)\ne\emptyset$, and let
$\alpha\colon G\to\Aut(A)$ be an action of a countable discrete group with equivariant property
(SI).
Then $\alpha$ is equivariantly $\mathcal{Z}$-stable if and only if, for every integer $m\ge2$,
there is a unital homomorphism
\[
 M_m\to(A^\omega\cap A')^{\alpha^\omega}.
\]
\end{thm}

Note that neither \autoref{thm:amenability-si} nor \autoref{thm:si-absorption} requires the
group to be amenable, so together they give the implication from \autoref{mainThmB} to
\autoref{mainThmA} when $T(A)\ne\emptyset$.

\subsection{Probability, joint moments, and spectral distributions}\label{sec:probability}
The construction in \autoref{sec:gaussian-gamma} uses selfadjoint elements whose moments
approach those of real Gaussian random variables.
We recall the basic terminology needed to read the main argument.
A self-contained introduction, including proofs of the joint moment estimates originating in
\cite{OAI}, is given in \autoref{sec:probability-appendix}.
That appendix develops all the probability notions recalled here, starting from scalar measures
and independence, and can be read before \autoref{sec:rows}.
The elements of our C*-algebras need not commute, so it is important to distinguish their
ordered products from products of scalar random variables.

\begin{notn}\label{not:probability}
Let $(\Omega,\mathcal F,\mathbb P)$ be a probability space.
For an integrable complex valued random variable $X$ on this space, we write
$\mathbb E(X)=\int_\Omega X\,d\mathbb P$ for its \emph{expectation}.
For a square integrable real random variable $X$, set
$\operatorname{Var}(X)=\mathbb E\bigl((X-\mathbb E(X))^2\bigr)$, and call this number its
\emph{variance}.
We say that $X$ is \emph{centered} if $\mathbb E(X)=0$.

The \emph{law} of a real valued random variable $X$ is the probability measure $\mu_X$ on $\R$ satisfying
\[
 \mu_X(V)=\mathbb P(\{\zeta\in\Omega\colon X(\zeta)\in V\})
 \quad\text{and}\quad
 \mathbb E(f(X))=\int_\R f\,d\mu_X
\]\nopagebreak
for all Borel sets $V\subseteq\R$ and all bounded Borel functions $f\colon\R\to\C$.
The \emph{moments} of $X$ are the numbers $\mathbb E(X^r)$, for integers $r\ge0$, whenever these integrals exist.
We say that its law is \emph{determined by its moments} if all these moments are finite and no
other probability measure on $\R$ has the same moments.

For a probability measure $\mu$ on $\R$, its \emph{Fourier transform} is the function
$\widehat{\mu}\colon\R\to\C$ given by
\[
 \widehat{\mu}(t)=\int_\R e^{itx}\,d\mu(x)
\]
for all $t\in\R$.
We use this sign convention throughout.
In probability theory this is also called the characteristic function of $\mu$, and we use
the term ``Fourier transform'' to distinguish it from the characteristic function of a set.
For the law of $X$, one has $\widehat{\mu_X}(t)=\mathbb E\bigl(e^{itX}\bigr)$ for all $t\in\R$.
\end{notn}

\Needspace{7\baselineskip}
\begin{df}\label{def:spectral-law}
Let $C$ be a unital C*-algebra, let $\rho$ be a state on $C$, and let $a\in C$ be selfadjoint.
The \emph{spectral distribution} of $a$ under $\rho$ is the probability measure $\mu_{a,\rho}$
on $\R$, supported on $\operatorname{sp}(a)$, such that
\[
 \rho(f(a))=\int_\R f\,d\mu_{a,\rho}
\]\nopagebreak
for all continuous functions $f\colon\R\to\C$.
\end{df}
The measure in \autoref{def:spectral-law} represents the restriction of $\rho$ to
$C^*(1,a)\cong C(\operatorname{sp}(a))$, extended to $\R$ by assigning measure zero to the
complement of $\operatorname{sp}(a)$.
In particular, its moments are $\rho(a^r)$ for integers $r\ge0$.
For a unitary $u\in C$, its spectral distribution is defined in the same way on $\T$, using
continuous functions on $\T$ and the moments $\rho(u^r)$ for $r\in\Z$.

\begin{rem}\label{rem:spectral-projections}
Suppose that a unitary $u\in B^\omega$ has the same spectral distribution $\mu$ under every $\tau\in T_\omega(B)$.
Let $V\subseteq\T$ be Borel.
Using regularity of $\mu$ and passing to a subsequence, choose continuous functions
$f_j\colon\T\to[0,1]$ converging to the characteristic function $\indicator_V$ both almost
everywhere and in $L^2(\mu)$.
Then
\[
 \uu{f_j(u)-f_k(u)}^2
 =\int_\T|f_j-f_k|^2\,d\mu
\]\nopagebreak
for all $j,k\in\N$.
Thus $(f_j(u))_{j\in\N}$ is Cauchy in the uniform $2$-norm, and completeness of the unit ball
gives a positive contraction $p\in B^\omega$ with $\lim_{j\to\infty}\uu{f_j(u)-p}=0$.
Since $\indicator_V^2=\indicator_V$ and the functions are bounded by one, dominated convergence gives
\[
 \lim_{j\to\infty}\uu{f_j(u)^2-f_j(u)}^2
 =\lim_{j\to\infty}\int_\T|f_j^2-f_j|^2\,d\mu=0.
\]
Continuity of multiplication on bounded sets now gives $p^2=p$, and $\tau(p)=\mu(V)$ for every
$\tau\in T_\omega(B)$.
Similarly, if $u$ commutes with a subset $S\subseteq B^\omega$, then so does $p$; if $u$ is
fixed by an action of a discrete group, then so is $p$.
Indeed, these properties hold for every $f_j(u)$ and pass to the uniform $2$-norm limit by \autoref{rem:ultranorm}.
This is the sense in which we apply characteristic functions of Borel sets to $u$ in
\autoref{sec:gamma}; no Borel functional calculus in $B^\omega$ is being assumed.
\end{rem}

\begin{df}\label{def:gaussian}
A \emph{standard real Gaussian variable} is a real random variable with probability density
$\frac{1}{\sqrt{2\pi}}e^{-\tfrac{t^2}{2}}$ on $\R$.
Let $J$ be a nonempty finite set.
A family $\xi=(\xi_g)_{g\in J}$ of real random variables is a \emph{centered real Gaussian
vector} if there are an integer $m\ge1$, independent standard real Gaussian variables
$\eta_1,\ldots,\eta_m$, and real numbers $a_{g,k}$, for $g\in J$ and $k=1,\ldots,m$, such that
$\xi$ has the same joint distribution as $\bigl(\sum_{k=1}^m a_{g,k}\eta_k\bigr)_{g\in J}$.
The coordinates $\xi_g$ are not required to be independent.
The \emph{covariance matrix} of $\xi$ is the matrix $R=(R_{g,h})_{g,h\in J}$, where
$R_{g,h}=\mathbb E(\xi_g\xi_h)$ for all $g,h\in J$.
\end{df}

\begin{rem}\label{rem:gaussian-covariance}
The covariance matrix of a centered real Gaussian vector is real symmetric and positive
semidefinite, and need not be invertible.
Conversely, let $R$ be a real symmetric positive semidefinite matrix indexed by $J$, and let
$\eta=(\eta_g)_{g\in J}$ have independent standard real Gaussian coordinates.
Then $R^{\tfrac{1}{2}}\eta$ is a centered real Gaussian vector with covariance $R$.
The joint law is determined by the covariance, and every fixed joint moment is a polynomial in
its entries; see \autoref{lem:app-wick}.
We give the construction and moment calculation in \autoref{app:gaussians}.
\end{rem}

We next specify what moment convergence means for elements of a C*-algebra, and recall the
estimates that will give it in our setting.

\begin{df}\label{def:ordered-moments}
Let $C$ be a unital C*-algebra with a tracial state $\tau$, let $J$ be a nonempty finite set,
and let $(a_g)_{g\in J}$ be selfadjoint elements of $C$.
An \emph{ordered word} in these elements is a product $a_{g_1}\cdots a_{g_d}$ with its
indicated order, and the corresponding \emph{joint $\tau$-moment} is
$\tau(a_{g_1}\cdots a_{g_d})$.
A \emph{noncommutative polynomial} is a finite linear combination of such words, including the
empty word $1$; for elements which are not selfadjoint, we also allow their adjoints in these
words.
For each $r\in\N$, let $C_r$ be a unital C*-algebra with a tracial state $\tau_r$, and let
$(a_g^{(r)})_{g\in J}$ be selfadjoint elements of $C_r$.
We say that their joint $\tau_r$-moments converge to the joint $\tau$-moments of
$(a_g)_{g\in J}$ as $r\to\infty$ if
\[
 \lim_{r\to\infty}\tau_r(a_{g_1}^{(r)}\cdots a_{g_d}^{(r)})
 =\tau(a_{g_1}\cdots a_{g_d})
\]
for every integer $d\ge1$ and all $g_1,\ldots,g_d\in J$.
By linearity, this is equivalent to convergence after evaluating each fixed noncommutative polynomial.
\end{df}

\begin{df}\label{def:row-algebras}
Let $C$ be a unital C*-algebra with a tracial state $\tau$.
A family $(C_r)_{r\in\N}$ of pairwise commuting unital C*-subalgebras of $C$ is
\emph{independent with respect to $\tau$} if
\[
 \tau(a_1\cdots a_s)=\tau(a_1)\cdots\tau(a_s)
\]
whenever $a_j\in C_{r_j}$ for $j=1,\ldots,s$ and the indices $r_1,\ldots,r_s$ are distinct.
Let $J$ be a nonempty finite set, and let $x_{r,g}\in C_r$ for all $r\in\N$ and all $g\in J$.
We call $(x_{r,g})_{g\in J}$ the $r$th \emph{row}, and also call $C_r$ a row algebra.
\end{df}
Elements in distinct rows commute, but elements within one row need not commute.

\begin{notn}\label{not:fluctuations}
Let $C$ be a unital C*-algebra and $\tau\in T(C)$. Let $(C_r)_{r\in\N}$ be pairwise commuting
unital C*-subalgebras of $C$ that are independent with respect to $\tau$.
Let $J$ be a nonempty finite set.
For $g\in J$ and $r\in\N$, let $x_{r,g}\in C_r$ be selfadjoint and satisfy $\tau(x_{r,g})=0$
and $\tau(x_{r,g}^2)=1$.
For $N\in\N$ and $g,h\in J$, set
\begin{equation}\label{eq:zk}
 z_{N,g}=\frac{1}{\sqrt N}\sum_{r=1}^Nx_{r,g}
 \quad\text{and}\quad
 k_{N;g,h}=\frac{1}{N}\sum_{r=1}^N\frac{x_{r,g}x_{r,h}+x_{r,h}x_{r,g}}{2}.
\end{equation}
We call $z_{N,g}$ a \emph{fluctuation variable} and $k_{N;g,h}$ a \emph{covariance variable}.
Set $R_N^\tau=(\tau(k_{N;g,h}))_{g,h\in J}$.
Traciality gives $\tau(k_{N;g,h})=\frac{1}{N}\sum_{r=1}^N\tau(x_{r,g}x_{r,h})$ for all
$g,h\in J$, so the entries of $R_N^\tau$ are real and symmetric.
If $t_g\in\R$ for $g\in J$, then
\begin{equation}\label{eq:empirical-covariance}
 \sum_{g,h\in J}t_gt_h\tau(k_{N;g,h})
 =\frac{1}{N}\sum_{r=1}^N\tau\Bigl(\bigl(\sum_{g\in J}t_gx_{r,g}\bigr)^2\Bigr)\ge0
\end{equation}\nopagebreak
Thus $R_N^\tau$ is a real symmetric positive semidefinite matrix with diagonal one, and is
therefore the covariance matrix of a centered real Gaussian vector by
\autoref{rem:gaussian-covariance}.
We call $k_{N;g,h}-\tau(k_{N;g,h})1$ a \emph{centered covariance variable}.
\end{notn}

The two estimates proved in \autoref{prop:moment-estimates} compare fluctuation moments with
Gaussian moments and bound words containing centered covariance factors.
Their constants depend only on the word length and the common norm bound for the row elements,
so they apply also when the trace varies.
In \autoref{thm:clt}, we verify independence for the row algebras arising from central copies
of $\mathcal{Z}$ and use these estimates to obtain joint moment convergence.

\section{Central weights and joint moment convergence}\label{sec:rows}

In this section, we construct a suitable commutative algebra and prove convergence of moments
for the construction of approximately invariant unitaries in \autoref{sec:gaussian-gamma}.
We first show that a tracially amenable action can be replaced by an auxiliary action admitting
a unital $G$-invariant central subalgebra on which the action is amenable.
We then construct commuting orbit algebras and apply the joint moment estimates from
\cite{OAI}, proved here in \autoref{sec:probability-appendix}.
Finally, we show how to include central coefficients in the moment calculations, allowing their
scalar values to vary with the extreme trace.

\subsection{A commutative amenable $G$-algebra in the relative commutant}\label{sec:reduction}

We begin by showing that the completely positive map obtained from tracial amenability can be
chosen to have commuting range.\footnote{Recall that a linear map $\Psi\colon C\to D$ is said
to have \emph{commuting range} if $[\Psi(c_1),\Psi(c_2)]=0$ for all $c_1,c_2\in C$.
This does not require $\Psi(c_1)\Psi(c_2)$ to belong to $\Psi(C)$.}
The algebra generated by its image will then provide the commutative algebra needed for the auxiliary action.
\begin{lma}\label{lem:commutative}
Let $A$ be a separable unital C*-algebra with $T(A)\ne\emptyset$, let $G$ be a countable
discrete group, and let $\alpha\colon G\to\Aut(A)$ be tracially amenable.
Then there is a separable, unital, commutative, $G$-invariant C*-subalgebra
$D\subseteq A^\omega\cap A'$ on which the induced action is amenable.
\end{lma}
	\begin{proof}
	Set $X=T(A)$, with the action induced by $\alpha$, and write $\gamma$ for the corresponding action on $C(X)$.
	The space $X$ is nonempty, compact, and metrizable, and the action on it is topologically
	amenable by \cite[Theorem 2.5]{GGKNV}.
	In particular, $G$ is exact, so the characterization of amenability using completely
	positive maps in the same theorem allows us to find an equivariant unital completely
	positive map
 $\Phi\colon C(X)\to A^\omega\cap A'$.
	By the Choi--Effros lifting theorem \cite{CE}, there are unital completely positive maps
	$\phi_j\colon C(X)\to A$, for $j\in\N$, representing $\Phi$.
	Centrality and equivariance of $\Phi$ give
	\begin{equation}\label{eq:phi-limits}
 \lim_{j\to\omega}\uu{[\phi_j(f),a]}=0
 \quad\text{and}\quad
 \lim_{j\to\omega}\uu{\alpha_h(\phi_j(f))-\phi_j(\gamma_h(f))}=0
	\end{equation}\nopagebreak
	for all $f\in C(X)$, all $a\in A$, and all $h\in G$.

	Choose increasing sequences $(F_n)_{n\in\N}$ and $(L_n)_{n\in\N}$ of finite subsets of
	$C(X)_1$ and $A_1$, respectively, whose unions are norm dense.
	Let $(K_n)_{n\in\N}$ be an increasing sequence of finite subsets of $G$ with union equal to $G$.
	For each $n\in\N$, we choose distinct indices $j(n,1),\ldots,j(n,n)$ such that
	\begin{equation}\label{eq:commutingchoices}
 \uu{[\phi_{j(n,r)}(f),a]}<\frac{1}{n}
 \quad\text{and}\quad
 \uu{\alpha_h(\phi_{j(n,r)}(f))-\phi_{j(n,r)}(\gamma_h(f))}<\frac{1}{n}
	\end{equation}\nopagebreak
	for all $f\in F_n$, all $a\in L_n$, all $h\in K_n$, and all $r=1,\ldots,n$, and
	\begin{equation}\label{eq:crosscommutators}
 \uu{[\phi_{j(n,s)}(f),\phi_{j(n,r)}(k)]}<\frac{1}{n}
	\end{equation}\nopagebreak
	for all $f,k\in F_n$ and all $1\le r<s\le n$.
	To see that this is possible, suppose that the first $s-1$ indices have been chosen.
	The elements $\phi_{j(n,r)}(k)$, for $r<s$ and $k\in F_n$, form a finite subset of $A$.
	By \eqref{eq:phi-limits}, the indices satisfying \eqref{eq:commutingchoices} and commuting
	within $\frac{1}{n}$ in uniform $2$-seminorm with this finite set form an $\omega$-large
	set.
	Since $\omega$ is free, we can choose $j(n,s)$ in this set and distinct from the preceding indices.

	For $n\in\N$, set $\psi_n=\frac{1}{n}\sum_{r=1}^n\phi_{j(n,r)}$, which is a unital
	completely positive map from $C(X)$ to $A$.
	Fix $f,k\in F_n$ and expand $[\psi_n(f),\psi_n(k)]$ as a double sum.
	The $n$ terms with $r=s$ have norm at most $2$.
	For $r>s$, the corresponding commutator is bounded by \eqref{eq:crosscommutators}; for
	$r<s$, apply the same estimate with $f$ and $k$ interchanged and use $[a,b]=-[b,a]$.
	Thus each of the remaining $n(n-1)$ terms has uniform $2$-norm less than $\frac{1}{n}$.
	Therefore
	\begin{align*}
 \uu{[\psi_n(f),\psi_n(k)]}
 &\le\frac1{n^2}\sum_{r,s=1}^n
       \uu{[\phi_{j(n,r)}(f),\phi_{j(n,s)}(k)]}\\
 &\mathrel{\overset{\mathclap{\eqref{eq:crosscommutators}}}{\le}}\quad\frac{2n}{n^2}+\frac{n(n-1)}{n^3}
 \le\frac3n.
	\end{align*}
	Similarly, the triangle inequality and \eqref{eq:commutingchoices} give
	\[
 \uu{[\psi_n(f),a]}<\frac{1}{n}
 \quad\text{and}\quad
 \uu{\alpha_h(\psi_n(f))-\psi_n(\gamma_h(f))}<\frac{1}{n}
	\]\nopagebreak
	for all $f\in F_n$, all $a\in L_n$, and all $h\in K_n$.
	Since these finite sets are increasing, the commutators and equivariance defects tend to
	zero for each fixed choice of elements in their unions, and contractivity and density give
	the same conclusion for arbitrary $f,k\in C(X)$, $a\in A$, and $h\in G$.

	Define $\Psi\colon C(X)\to A^\omega\cap A'$ by $\Psi(f)=[(\psi_n(f))_{n\in\N}]$ for all $f\in C(X)$.
	The preceding estimates show that $\Psi$ is equivariant and has commuting range.
	Moreover, $\Psi$ is unital and completely positive because each $\psi_n$ is.
	Since $\Psi$ preserves adjoints, its range generates a unital separable commutative
	C*-algebra $D=C^*(\Psi(C(X)))$.
	Moreover, equivariance gives $\alpha_h^\omega(\Psi(C(X)))=\Psi(C(X))$ for every $h\in G$, so $D$ is $G$-invariant.

	It remains to verify amenability of the induced action on $D$.
	For a real valued $f\in C(X)$, positivity of $\Psi$ and commutativity of $D$ give
	\begin{equation}\label{eq:positive-absolute}
 |\Psi(f)|\le\Psi(|f|).
	\end{equation}
	Let $K\subseteq G$ be finite, and let $\varepsilon>0$.
	Choose amenability weights $(m_g)_{g\in G}$ for $C(X)$ satisfying \eqref{eq:weights} for
	$K$ and $\varepsilon$, and set $d_g=\Psi(m_g)$ for all $g\in G$.
	The family $(d_g)_{g\in G}$ has finite support, consists of positive elements, and sums to
	one because $\Psi$ is positive and unital.
	For every $h\in K$, all the differences $\gamma_h(m_{h^{-1}g})-m_g$ are selfadjoint.
	Using equivariance, we obtain
	\begin{align*}
 \sum_{g\in G}|\alpha_h^\omega(d_{h^{-1}g})-d_g|
 ={}&\sum_{g\in G}|\Psi(\gamma_h(m_{h^{-1}g})-m_g)|\\
 \stackrel{\mathclap{\eqref{eq:positive-absolute}}}{\le}{}&\Psi\Bigl(\sum_{g\in G}|\gamma_h(m_{h^{-1}g})-m_g|\Bigr).
	\end{align*}
	Taking norms and using contractivity of $\Psi$, we conclude that $(d_g)_{g\in G}$ satisfies
	\eqref{eq:weights} for $K$ and $\varepsilon$.
	Since $K$ and $\varepsilon$ were arbitrary, the induced action on $D$ is amenable.
\end{proof}

The preceding lemma gives an amenable commutative algebra in the relative commutant, but we
need the weights to be central in the algebra on which we work.
We achieve this by passing to a tensor product, while retaining a homomorphism back to the original ultrapower.

For a unital C*-algebra $A$ with $T(A)\ne\emptyset$, denote by $\Delta_A\colon A\to A^\omega$
the homomorphism given by constant sequences.

\begin{prop}\label{prop:auxiliary}
Let $A$ be a separable, unital, $\mathcal{Z}$-stable C*-algebra with $T(A)\ne\emptyset$.
Let $\alpha\colon G\to\Aut(A)$ be a tracially amenable action of a countable discrete group.
Then there are a separable, unital, $\mathcal{Z}$-stable C*-algebra $B$ with
$T(B)\ne\emptyset$, an action $\beta\colon G\to\Aut(B)$, a unital $G$-invariant subalgebra
$E\subseteq Z(B)$ on which the action is amenable, and equivariant unital homomorphisms
$\iota\colon A\to B$ and $\theta\colon B\to A^\omega$ making the following diagram commute.
\[
\begin{tikzcd}[column sep=large,row sep=large]
 A \arrow[rr,"\Delta_A"] \arrow[dr,"\iota"'] && A^\omega \\
 & B \arrow[ur,"\theta"'] &
\end{tikzcd}
\]
\end{prop}
	\begin{proof}
	Choose $D$ as in \autoref{lem:commutative}, and denote by $\gamma$ the action $\alpha^\omega|_D$.
	Set $B=A\otimes D$, set $\beta=\alpha\otimes\gamma$, and set $E=1\otimes D$.
	The algebra $B$ is separable, unital, and $\mathcal{Z}$-stable. It has a tracial state: for
	$\tau\in T(A)$ and any character $\chi$ of the commutative unital algebra $D$, the tensor
	product $\tau\otimes\chi$ belongs to $T(B)$.
	Moreover, commutativity of $D$ gives $E\subseteq Z(B)$, and the restriction of $\beta$ to
	$E$ is amenable because it is conjugate to $\gamma$.

	Since $D$ is commutative, it is nuclear, and hence $A\otimes_{\max}D=A\otimes D$. The
	commuting homomorphisms $\Delta_A\colon A\to A^\omega$ and $D\hookrightarrow A^\omega$
	therefore define a unital homomorphism $\theta\colon B\to A^\omega$ satisfying
	$\theta(a\otimes d)=\Delta_A(a)d$ for all $a\in A$ and all $d\in D$.
	For $g\in G$, $a\in A$, and $d\in D$, we have
	\[
 \theta(\beta_g(a\otimes d))
 =\Delta_A(\alpha_g(a))\alpha_g^\omega(d)
 =\alpha_g^\omega(\theta(a\otimes d)),
	\]
	so $\theta$ is equivariant.
	Finally, define $\iota\colon A\to B$ by $\iota(a)=a\otimes1$ for all $a\in A$.
	Then $\iota$ is an equivariant unital embedding, and $(\theta\circ\iota)(a)=\Delta_A(a)$ for all $a\in A$.
\end{proof}

\subsection{Commuting rows and the moment theorem}

We next recall the construction of commuting orbit algebras and the moment estimates from \cite{OAI}.
Since we will use these estimates after adjoining the central amenability weights, we also
check how the extreme traces change when the central algebra is taken into account.

\begin{prop}\label{prop:rows}
Let $B$ be a separable, unital, $\mathcal{Z}$-stable C*-algebra with $T(B)\ne\emptyset$, and
let $\beta\colon G\to\Aut(B)$ be an action of a countable discrete group.
Let $L\subseteq B^\omega$ be a norm separable $G$-invariant C*-subalgebra.
Then there are unital homomorphisms $\iota_r\colon\mathcal{Z}\to B^\omega$, for $r\in\N$, such that the algebras
\[
 C_r=C^*(\{\beta_g^\omega(\iota_r(\mathcal{Z}))\colon g\in G\})
\]
are $G$-invariant, commute with $L$, and commute pairwise.
\end{prop}
	\begin{proof}
	Use $\mathcal{Z}$-stability of $B$ and \cite[Theorem 2.2]{TW} to find a unital homomorphism
	$\Xi\colon\mathcal{Z}\to B_\omega\cap B'$.
	Use the Choi--Effros theorem to choose a unital completely positive lift
	$\Psi\colon\mathcal{Z}\to\ell^\infty(\N,B)$.
	Write $\Psi=(\psi_n)_{n\in\N}$ and note that the maps $\psi_n\colon\mathcal{Z}\to B$ are
	asymptotically multiplicative and central in norm along $\omega$.
\vspace{.2cm}

	\textbf{Claim:} \emph{For every norm separable C*-subalgebra $L_0\subseteq B^\omega$, there is a unital homomorphism $\iota\colon\mathcal{Z}\to B^\omega\cap L_0'$.} Choose a norm dense sequence $(z_r)_{r\in\N}$ in the unit ball of $\mathcal{Z}$, and a norm dense sequence $(a_r)_{r\in\N}$ in the unit ball of $L_0$.
	For each $r\in\N$, fix a bounded sequence $(a_{r,k})_{k\in\N}$ representing $a_r$.
	For each $k\in\N$, the elements $a_{1,k},\ldots,a_{k,k}$ are fixed elements of $B$.
	By asymptotic centrality and multiplicativity, the indices satisfying the following
	inequalities form an $\omega$-large set, so we can choose $n(k)\ge k$ such that
	\begin{equation}\label{eq:row-reindex}
 \|[\psi_{n(k)}(z_r),a_{s,k}]\|<\frac{1}{k}
 \quad\text{and}\quad
 \|\psi_{n(k)}(z_rz_s)-\psi_{n(k)}(z_r)\psi_{n(k)}(z_s)\|<\frac{1}{k}
	\end{equation}\nopagebreak
	for all $r,s=1,\ldots,k$.

	Define $\iota\colon\mathcal{Z}\to B^\omega$ by $\iota(z)=[(\psi_{n(k)}(z))_{k\in\N}]$ for all $z\in\mathcal{Z}$.
	The map is unital and completely positive, and \eqref{eq:row-reindex} gives
	$\iota(z_rz_s)=\iota(z_r)\iota(z_s)$ for all $r,s\in\N$, since norm convergence implies
	uniform $2$-norm convergence.
	By contractivity and density, $\iota$ is multiplicative on all of $\mathcal{Z}$, and hence is a homomorphism.
	The commutator estimate in \eqref{eq:row-reindex}, together with density, shows that its
	range commutes with $L_0$.
	\vspace{.2cm}

	We apply the claim successively to the $G$-invariant algebra generated by $L$ and the
	previously constructed rows, which is norm separable because $G$ is countable.
	To check commutation of the whole new row, let $a$ belong to this $G$-invariant algebra,
	let $z\in\mathcal{Z}$, and let $g\in G$.
	Then
	\[
 [\beta_g^\omega(\iota(z)),a]
 =\beta_g^\omega([\iota(z),\beta_{g^{-1}}^\omega(a)])=0.
	\]
	Thus every translate of the new copy commutes with $L$ and all the earlier rows.
	The algebra generated by these translates is $G$-invariant by construction.
	Induction gives the required family $(C_r)_{r\in\N}$.
\end{proof}

We record the following observation to obtain row elements with the same normalization under
every trace, as required by \autoref{prop:moment-estimates}.
\begin{rem}\label{rem:row-elements}
Adopt the assumptions and notation from \autoref{prop:rows}.
Since $\tau_{\mathcal{Z}}$ is faithful and $\mathcal{Z}\ne\C$, we can choose a selfadjoint
element $b\in\mathcal{Z}$ with $\tau_{\mathcal{Z}}(b)=0$ and $\tau_{\mathcal{Z}}(b^2)=1$.
Set $x_{r,g}=\beta_g^\omega(\iota_r(b))$ for all $r\in\N$ and all $g\in G$, and set
$C_0=C^*(\bigcup_{r\in\N}C_r)$ and $M=\|b\|$.
Every tracial state on $C_0$ restricts to the unique trace on each translated copy of $\mathcal{Z}$.
Hence
\begin{equation}\label{eq:rowmoments}
 \|x_{r,g}\|\le M,\quad \tau(x_{r,g})=0,
 \quad\text{and}\quad
 \tau(x_{r,g}^2)=1
\end{equation}\nopagebreak
for all $r\in\N$, $g\in G$, and $\tau\in T(C_0)$.
We also have $\beta_h^\omega(x_{r,g})=x_{r,hg}$ for all $h,g\in G$.
\end{rem}

The following is the main result of this section, which we will informally call the moment theorem.
It extends \cite[Proposition 3.3]{OAI} to include a central subalgebra, with no assumptions on
simplicity or on a group action.
We state the central coefficients and their limits explicitly, since this is the form used in
\autoref{sec:gaussian}.

\begin{thm}\label{thm:clt}
Let $C$ be a unital C*-algebra, and let $E\subseteq Z(C)$ be a unital C*-subalgebra.
Let $(C_r)_{r\in\N}$ be pairwise commuting unital C*-subalgebras of $C$.
Set $C_0=C^*(\bigcup_{r\in\N}C_r)$, and suppose that $C=C^*(E,C_0)$.
Let $J$ be a nonempty finite set, let $M>0$, and let $x_{r,g}\in C_r$ be selfadjoint elements satisfying
\[
 \|x_{r,g}\|\le M,\quad \tau(x_{r,g})=0,
 \quad\text{and}\quad \tau(x_{r,g}^2)=1
\]
for all $r\in\N$, all $g\in J$, and all $\tau\in T(C_0)$.
Define $z_{N,g}$ and $k_{N;g,h}$ by \eqref{eq:zk}.
Let $(N_j)_{j\in\N}$ be a sequence of positive integers tending to infinity, let
$(\tau_j)_{j\in\N}$ be a sequence in $\partial_eT(C)$, let $s\in\N$, and let
$e_1,\ldots,e_s\in E$.
Suppose that all the following limits exist, and set
\[
 R_{g,h}=\lim_{j\to\infty}\tau_j(k_{N_j;g,h})
 \quad\text{and}\quad
 \lambda_k=\lim_{j\to\infty}\tau_j(e_k)
\]
for all $g,h\in J$ and all $k=1,\ldots,s$.
Set $R=(R_{g,h})_{g,h\in J}$.
Then $R$ is real symmetric and positive semidefinite with diagonal one.
Let $\xi=(\xi_g)_{g\in J}$ be a centered real Gaussian vector with covariance $R$, and set
$\mathbf z_N=(z_{N,g})_{g\in J}$ and $\mathbf k_N=(k_{N;g,h})_{g,h\in J}$ for all $N\in\N$.
When evaluating a noncommutative polynomial, we regard a matrix as the family of its entries.
For every noncommutative polynomial $P$ in these variables and their adjoints, one has
\[
 \lim_{j\to\infty}\tau_j\bigl(P(\mathbf z_{N_j},\mathbf k_{N_j},e_1,\ldots,e_s)\bigr)
 =\mathbb E\bigl(P(\xi,R,\lambda_1,\ldots,\lambda_s)\bigr).
\]
For any sequences $(N_j)_{j\in\N}$ and $(\tau_j)_{j\in\N}$ as above, there is a subsequence
along which the covariance entries and the values of the fixed coefficients $e_1,\ldots,e_s$
all converge.
\end{thm}
	\begin{proof}
	We first explain the role of centrality.
	Fix $\tau\in\partial_eT(C)$, and set $M_0=\pi_\tau(C)''$ and $M_r=\pi_\tau(C_r)''$ for all $r\in\N$.
	Since $M_0$ is a factor and $E$ is central, $\pi_\tau(e)=\tau(e)1$ for all $e\in E$.
	Since $C=C^*(E,C_0)$, it follows that $\pi_\tau(C_0)''=M_0$, so $\tau|_{C_0}$ is an extreme trace on $C_0$.
	Denote by $\widehat\tau$ the faithful normal tracial state on $M_0$ extending $\tau$.
	The algebras $M_r$ commute pairwise and generate $M_0$, so every element of $Z(M_r)$ is central in $M_0$.
	Thus each $M_r$ is a factor.

	There is a unique $\widehat\tau$-preserving conditional expectation $E_r\colon M_0\to M_r$;
	see \cite[Theorem IX.4.2]{Takesaki}.
	If $y$ belongs to the von Neumann algebra generated by the other rows, bimodularity gives
	\[
 aE_r(y)=E_r(ay)=E_r(ya)=E_r(y)a
	\]
	for every $a\in M_r$.
	Hence $E_r(y)\in Z(M_r)=\C1$, and preservation of the trace gives $E_r(y)=\widehat\tau(y)1$.
	It follows that
	\[
 \widehat\tau(ay)=\widehat\tau(E_r(ay))
 =\widehat\tau(aE_r(y))
 =\widehat\tau(a)\widehat\tau(y)
	\]
	for every $a\in M_r$ and every such $y$.
	Applying this identity successively, we obtain
	\[
 \tau(a_1\cdots a_t)=\tau(a_1)\cdots\tau(a_t)
	\]
	whenever $a_i\in C_{r_i}$ for $i=1,\ldots,t$ and the row indices $r_1,\ldots,r_t$ are distinct.
	This proves the independence required by \autoref{prop:moment-estimates}.

	For an extreme trace $\tau$ on $C$, set $R_N^\tau=(\tau(k_{N;g,h}))_{g,h\in J}$.
	As in \autoref{not:fluctuations}, these matrices are real symmetric and positive
	semidefinite with diagonal one, so the limit matrix $R$ has the same properties; see
	\autoref{rem:gaussian-covariance}.
	Independence and the bounds on $x_{r,g}$ verify the hypotheses of
	\autoref{prop:moment-estimates}, with constants which do not depend on $\tau$.

	We next prove convergence of joint moments for the fluctuation and covariance variables.
	Fix an ordered monomial $P$ containing $n$ fluctuation variables and $m$ covariance variables.
	For $N\in\N$ and $g,h\in J$, set $c_{N;g,h}=k_{N;g,h}-\tau(k_{N;g,h})1$, which we call a
	centered covariance variable.
	In $P$, replace each covariance variable by
	\[
 k_{N;g,h}=\tau(k_{N;g,h})1+c_{N;g,h},
	\]
	and expand the product, preserving the order of its factors.
	There are $2^m$ terms, according to the choice of a scalar or a centered covariance variable at each position.
	The scalar coefficients have absolute value at most one by Cauchy--Schwarz.
	By \eqref{eq:centered-error}, the absolute value of the trace of a term with $p\ge1$
	centered covariance factors is at most $K_{n,p,M}/N^{p/2}$, which tends to zero uniformly
	in $\tau$.
	There is only one term left, namely the one in which all covariance factors are replaced by their scalar traces.
	When $n\ge1$, its remaining fluctuation moment differs from the corresponding Gaussian
	moment by at most $K_{n,M}/\sqrt N$, by \eqref{eq:moment-error}; when $n=0$, the remaining
	factor is $1$.
	Applying this argument with $N=N_j$ and $\tau=\tau_j$, we find that the scalar coefficients
	converge to the entries of $R$, while \autoref{lem:app-wick} gives convergence of the
	Gaussian moments.
	This proves the assertion for $P$, and linearity gives it for every noncommutative
	polynomial in these two families.

	Finally, in the $j$th tracial representation each $e_k$ is the scalar $\tau_j(e_k)$.
	Thus an ordered word which also contains the $e_k$ and their adjoints has trace equal to a
	row moment multiplied by the corresponding product of these scalars and their conjugates.
	Their limits are $\lambda_k$ and $\overline{\lambda_k}$, respectively, so the convergence
	just proved gives the displayed formula for every polynomial involving the central
	coefficients as well.
	The covariance entries lie in $[-1,1]$, and $|\tau_j(e_k)|\le\|e_k\|$ for all $j\in\N$ and $k=1,\ldots,s$.
	Compactness of this finite collection of bounded scalar parameters gives the last assertion.
\end{proof}

\section{Gaussian quadratic forms and local equivariant property Gamma}\label{sec:gaussian-gamma}

In this section, we prove that an action on a separable, unital, $\mathcal{Z}$-stable C*-algebra
with nonempty tracial state space has local equivariant property Gamma whenever it admits a
unital $G$-invariant central subalgebra on which the action is amenable; see
\autoref{thm:gamma}.
Using the moment convergence from \autoref{sec:rows}, we first construct approximately
invariant unitaries whose nonzero moments are uniformly small over the trace space.
We then pass to projections that divide every limit trace equally and commute with a prescribed
separable subset of the ultrapower.

\subsection{Gaussian quadratic forms with central weights}\label{sec:gaussian}

Our next goal is to use \autoref{thm:clt} to construct approximately invariant unitaries whose
nonzero moments are small.
Since a weighted sum of the Gaussian variables can have arbitrarily small variance, we use
quadratic forms and apply a square root to obtain a lower bound that also holds when the
covariance matrix is singular.

For operators on a finite dimensional real Euclidean space $V$, we use the unnormalized matrix
trace $\Tr$, the Hilbert--Schmidt norm $\HS T=\Tr(T^*T)^{\tfrac{1}{2}}$, and the trace norm
$\|T\|_1=\Tr(|T|)$.
A standard Gaussian vector on $V$ has independent standard real Gaussian coordinates in an orthonormal basis.
For such a vector $\eta$ and a selfadjoint $T$, set $Q_T=\langle\eta,T\eta\rangle-\Tr(T)$.
The estimates for these quadratic forms, the passage from moments to continuous functional
calculus, and the perturbation estimate for exponentials are recalled and proved in
\autoref{sec:probability-appendix}; the assertions used below are collected in
\autoref{prop:gaussianinput}.

The next proposition extends \cite[Proposition 4.1]{OAI} to central amenability weights.
Their values at characters vary with the extreme trace, so the estimates must be uniform in
both the covariance matrix and the scalar weights.

\begin{prop}\label{prop:unitaries}
Let $G$ be a countable discrete group, let $C$ be a unital C*-algebra with $T(C)\ne\emptyset$,
and let $\gamma\colon G\to\Aut(C)$ be an action.
Suppose that $C$ is generated by a unital $G$-invariant subalgebra $E\subseteq Z(C)$ and
pairwise commuting unital $G$-invariant subalgebras $C_r$, for $r\in\N$.
Set $C_0=C^*(\bigcup_{r\in\N}C_r)$.
Assume that the restriction of $\gamma$ to $E$ is amenable, and that there are selfadjoint
elements $x_{r,g}\in C_r$ and a constant $M>0$ satisfying $\gamma_h(x_{r,g})=x_{r,hg}$ and
$\|x_{r,g}\|\le M$ for all $r\in\N$ and $g,h\in G$, as well as
\[
 \tau(x_{r,g})=0
 \quad\text{and}\quad\tau(x_{r,g}^2)=1
\]\nopagebreak
for all $r\in\N$, all $g\in G$, and all $\tau\in T(C_0)$.
For every finite subset $K\subseteq G$, every integer $m\ge1$, and every $\varepsilon>0$, there
is a unitary $u\in C$ such that
\begin{equation}\label{eq:unitarytarget}
 \sup_{\tau\in T(C)}|\tau(u^d)|<\varepsilon
 \quad\text{and}\quad
 \sup_{\tau\in T(C)}\|\gamma_h(u)-u\|_{2,\tau}<\varepsilon
\end{equation}\nopagebreak
for all $d\in\Z$ with $1\le|d|\le m$ and all $h\in K$.
Moreover, $u$ can be chosen in the form $e^{ia}$ for a selfadjoint $a\in C$.
\end{prop}
	\begin{proof}
	We will construct selfadjoint elements whose joint tracial moments with their translates
	approach those of Gaussian quadratic forms.
	These quadratic forms have variance bounded away from zero and change little under the
	elements of $K$, uniformly over the traces.
	Exponentiating the selfadjoint elements will then give the desired unitary.
	Without loss of generality, we may assume that $0<\varepsilon<1$.
	Choose $t>0$ such that $(1+t^2)^{-\tfrac{1}{4}}<\frac{\varepsilon}{4}$.
	Next choose $\delta>0$ and $0<\varepsilon_0<\frac{1}{2}$ satisfying
	\begin{equation}\label{eq:choices}
 t\sqrt2(\sqrt\delta+2\varepsilon_0)<\frac{\varepsilon}{4}.
	\end{equation}
	Choose amenability weights $(d_g)_{g\in G}$ in $E$ satisfying \eqref{eq:weights} with
	$\delta$ in place of $\varepsilon$.
	Set $L=\{g\in G\colon d_g\ne0\}$ and $n=|L|$.
	Set $J=L\cup\bigcup_{h\in K}hL$.
	For $h\in K$ and $g\in G$, set $d_g^{(h)}=\gamma_h(d_{h^{-1}g})$.
	Both $(d_g)_{g\in G}$ and $(d_g^{(h)})_{g\in G}$ have sum one, and
	\begin{equation}\label{eq:weightdefect}
 \Bigl\|\sum_{g\in G}|d_g^{(h)}-d_g|\Bigr\|<\delta
	\end{equation}\nopagebreak
	for all $h\in K$.

	We first establish estimates for a scalar Gaussian model.
	Let $R=(R_{g,r})_{g,r\in J}$ be a real positive semidefinite matrix with diagonal one, and
	let $\chi$ be a character of $E$.
	Choose a finite dimensional real Euclidean space $V$ and unit vectors $v_g\in V$, for
	$g\in J$, which span $V$ and satisfy $\langle v_g,v_r\rangle=R_{g,r}$ for all $g,r\in J$.
	Such a family is called a family of Gram vectors for $R$; its existence is explained in \autoref{app:gaussians}.
	For $r,g\in J$, define the rank one operator $\theta_{r,g}\colon V\to V$ by
	$\theta_{r,g}(w)=\langle w,v_r\rangle v_g$ for all $w\in V$.
	For $h\in K$, set
	\[
 P=\sum_{g\in L}\chi(d_g)\theta_{g,g}
 \quad\text{and}\quad
 P_h=\sum_{g\in hL}\chi(d_g^{(h)})\theta_{g,g}.
	\]
	The coefficients are nonnegative and the vectors have norm one. Since each rank one
	operator $\theta_{g,g}$ has trace one and both families of weights sum to one, $P$ and
	$P_h$ are positive with $\Tr(P)=\Tr(P_h)=1$. The spectra of $P$ and $P_h$ therefore lie in
	$[0,1]$, and their ranks are at most $n$ because each is a sum of at most $n$ rank one
	operators.
	Since each rank one matrix $\theta_{g,g}$ has trace norm one, we have
	\[
 \|P-P_h\|_1\le\sum_{g\in J}|\chi(d_g)-\chi(d_g^{(h)})|
 \stackrel{\eqref{eq:weightdefect}}{<}\delta.
	\]
	\autoref{prop:gaussianinput}\ref{item:sqrt} now implies that
	\begin{equation}\label{eq:Pbounds}
 \Bigl\|P^{\frac{1}{2}}\Bigr\|_{\mathrm{HS}}=1
 \quad\text{and}\quad \Bigl\|P^{\frac{1}{2}}-P_h^{\frac{1}{2}}\Bigr\|_{\mathrm{HS}}<\sqrt\delta.
	\end{equation}

	Choose a real polynomial $p$ with $p(0)=0$ such that
	\begin{equation}\label{eq:polynomialapprox}
 \sup_{0\le x\le1}|p(x)-\sqrt x|<\frac{\varepsilon_0}{\sqrt n}.
	\end{equation}
	Write $p(x)=\sum_{r=1}^Dc_rx^r$.
	Let $S$ be either $P$ or $P_h$.
	Since $p(0)=0$, the operator $p(S)-S^{\tfrac{1}{2}}$ vanishes on the kernel of $S$ and has rank at most $n$.
	Diagonalizing $S$, we obtain
	\[
 \Bigl\|p(S)-S^{\frac{1}{2}}\Bigr\|_{\mathrm{HS}}^2
 \le n\sup_{0\le x\le1}|p(x)-\sqrt x|^2
 \stackrel{\eqref{eq:polynomialapprox}}{<}\varepsilon_0^2.
	\]
	Set $T=p(P)$ and $T_h=p(P_h)$.
	The triangle inequality gives
	\begin{equation}\label{eq:Tbounds}
	\begin{split}
 \HS T&\ge\Bigl\|P^{\frac{1}{2}}\Bigr\|_{\mathrm{HS}}-\Bigl\|T-P^{\frac{1}{2}}\Bigr\|_{\mathrm{HS}}
 \stackrel{\eqref{eq:Pbounds}}{>}1-\varepsilon_0>\frac{1}{2},\\
 \HS{T-T_h}&\le\Bigl\|T-P^{\frac{1}{2}}\Bigr\|_{\mathrm{HS}}+\Bigl\|P^{\frac{1}{2}}-P_h^{\frac{1}{2}}\Bigr\|_{\mathrm{HS}}
                  +\Bigl\|P_h^{\frac{1}{2}}-T_h\Bigr\|_{\mathrm{HS}}\\
 &\stackrel{\eqref{eq:Pbounds}}{<}\sqrt\delta+2\varepsilon_0.
	\end{split}
	\end{equation}
	Note that the bounds in \eqref{eq:Tbounds} do not depend on the character or the covariance matrix.
	On the other hand, the polynomial $p$ depends on $n$ and $\varepsilon_0$, and will remain
	fixed as the trace and $N$ vary.

	We next construct elements of $C$ whose moments converge to those of these Gaussian quadratic forms.
	Define $z_{N,g}$ and $k_{N;g,r}$ by \eqref{eq:zk}, for all $g,r\in G$.
	In the following formula, an empty product of covariance factors is understood to be the unit.
	For each $N\in\N$, set
	\begin{equation}\label{eq:empirical}
	\begin{split}
 \widetilde a_N={}&\sum_{r=1}^Dc_r\sum_{g_1,\ldots,g_r\in L}
 d_{g_1}\cdots d_{g_r}
 k_{N;g_1,g_2}\cdots k_{N;g_{r-1},g_r}\\
 &\hspace{12mm}\cdot\bigl(z_{N,g_1}z_{N,g_r}-k_{N;g_r,g_1}\bigr).
	\end{split}
	\end{equation}
	The weights $d_{g_1},\ldots,d_{g_r}$ are central, but the fluctuation and covariance variables need not commute.
	Set $a_N=\frac{1}{2}(\widetilde a_N+\widetilde a_N^*)$, which is selfadjoint.
	The sequence $(a_N)_{N\in\N}$ need not be bounded in norm, so we will use the moment
	convergence statement in \autoref{prop:gaussianinput}\ref{item:spectral}.

	We now identify the Gaussian quadratic form whose moments are the limits of the tracial
	moments of this expression.
	Let $(N_j)_{j\in\N}$ be a sequence tending to infinity, and let $(\tau_j)_{j\in\N}$ be a
	sequence in $\partial_eT(C)$.
	Let $E_0$ be the unital commutative algebra generated by the finitely many weights $d_g$
	and $d_g^{(h)}$ occurring above.
	Since the character space of $E_0$ is compact and metrizable and the covariance entries are
	bounded, after passing to a subsequence we may assume that each sequence
	$(\tau_j(k_{N_j;g,r}))_{j\in\N}$ converges and that the characters $\tau_j|_{E_0}$ converge
	to a character $\chi$ of $E_0$.
	Denote the limiting covariance by $R$.
	All matrices and estimates above only use elements of $E_0$, so they apply to this character.

	Choose $V$ and Gram vectors $(v_g)_{g\in J}$ for this limiting matrix $R$ as above, and
	form the corresponding operators $P$, $P_h$, $T$, and $T_h$ using $\chi$.
	Let $\eta$ be a standard Gaussian vector on $V$, and set $\xi_g=\langle\eta,v_g\rangle$ for all $g\in J$.
	Expanding the power of $P$ gives
	\[
 P^r=\sum_{g_1,\ldots,g_r\in L}
 \chi(d_{g_1})\cdots\chi(d_{g_r})
 R_{g_1,g_2}\cdots R_{g_{r-1},g_r}\,\theta_{g_r,g_1}.
	\]
	For the final rank one operator in this expression, we have
	\[
 \langle \eta,\theta_{g_r,g_1}\eta\rangle=\xi_{g_1}\xi_{g_r}
 \quad\text{and}\quad
 \Tr(\theta_{g_r,g_1})=R_{g_r,g_1}.
	\]
	By \autoref{thm:clt}, the $r$th summand of \eqref{eq:empirical} has joint moments converging to those of
	$c_r\bigl(\langle \eta,P^r\eta\rangle-\Tr(P^r)\bigr)$.
	Adding over $r$ gives $\langle \eta,p(P)\eta\rangle-\Tr(p(P))=Q_T$.
	Every power of $\widetilde a_{N_j}$, or ordered product involving its adjoint, expands into
	finitely many words covered by \autoref{thm:clt}.
	The limiting expression $Q_T$ is real, so replacing $\widetilde a_{N_j}$ by its selfadjoint
	part does not change that limit.
	Therefore all moments of $a_{N_j}$ converge to the corresponding moments of $Q_T$.
	We next check the translated expression.
	Equivariance of the row variables gives
	\[
 \gamma_h(z_{N,g})=z_{N,hg}
 \quad\text{and}\quad
 \gamma_h(k_{N;g,r})=k_{N;hg,hr}
	\]\nopagebreak
	for all $g,r,h\in G$ and all $N\in\N$.
	Apply $\gamma_h$ to \eqref{eq:empirical}, and replace each summation index $g_j\in L$ by $hg_j\in hL$.
	The coefficient $\gamma_h(d_{g_j})$ is $d_{hg_j}^{(h)}$.
	Thus the same computation, with $hL$ and $d_g^{(h)}$ in place of $L$ and $d_g$, gives $Q_{T_h}$ as the limit.
	The joint ordered moments of $a_{N_j}$ and $\gamma_h(a_{N_j})$ therefore converge to those
	of $Q_T$ and $Q_{T_h}$, with the same $\eta$.
	In particular, $Q_{T_h}-Q_T=Q_{T_h-T}$, so we can apply \eqref{eq:gaussiandifference}.

	It remains to deduce uniform estimates for the unitaries.
	Set $u_N=e^{ita_N}$.
	Using \autoref{prop:gaussianinput}\ref{item:spectral} and the bound $\HS T>\frac{1}{2}$
	from \eqref{eq:Tbounds}, we get
	\[
 \lim_{j\to\infty}|\tau_j(u_{N_j}^d)|
 =|\mathbb{E}\bigl(e^{idtQ_T}\bigr)|
 \stackrel{\eqref{eq:gaussianest}}{\le}(1+t^2)^{-\frac{1}{4}}<\frac{\varepsilon}{4}
	\]\nopagebreak
	for all $d\in\Z$ with $1\le|d|\le m$.
	To obtain the invariance estimate, expand the square of the selfadjoint difference.
	The joint moment convergence just proved gives
	\begin{equation}\label{eq:quadraticdifference-limit}
 \lim_{j\to\infty}\bigl\|\gamma_h(a_{N_j})-a_{N_j}\bigr\|_{2,\tau_j}^2
 =\mathbb E\bigl((Q_{T_h}-Q_T)^2\bigr)
 =\|Q_{T_h}-Q_T\|_{L^2}^2
	\end{equation}\nopagebreak
	for every $h\in K$.
	Since $\gamma_h(u_N)=e^{it\gamma_h(a_N)}$, we obtain
	\begin{align*}
 \limsup_{j\to\infty}\bigl\|\gamma_h(u_{N_j})-u_{N_j}\bigr\|_{2,\tau_j}
 &\stackrel{\eqref{eq:duhamel}}{\le}
 t\lim_{j\to\infty}\bigl\|\gamma_h(a_{N_j})-a_{N_j}\bigr\|_{2,\tau_j}\\
 &\stackrel{\eqref{eq:quadraticdifference-limit}}{=}
 t\|Q_{T_h}-Q_T\|_{L^2}
 \stackrel{\eqref{eq:gaussiandifference}}{=}t\sqrt2\HS{T_h-T}\\
 &<t\sqrt2(\sqrt\delta+2\varepsilon_0)
 \stackrel{\eqref{eq:choices}}{<}\frac{\varepsilon}{4}.
	\end{align*}

	We claim that the desired estimates, with $\frac{\varepsilon}{2}$ in place of
	$\varepsilon$, hold uniformly over the extreme traces for all sufficiently large $N$.
	Arguing by contradiction, suppose otherwise and find a sequence $(N_j)_{j\in\N}$ tending to
	infinity and a sequence $(\tau_j)_{j\in\N}$ of extreme traces for which one of these
	estimates fails.
	Since there are only finitely many inequalities, we may pass to a subsequence on which one fixed inequality fails.
	Passing to a further subsequence gives the covariance and weight limits used above,
	contradicting the estimates just proved.

	For a fixed $N$, the evaluations $\tau(u_N^d)$ and
	$\tau\bigl((\gamma_h(u_N)-u_N)^*(\gamma_h(u_N)-u_N)\bigr)$, viewed as functions of
	$\tau\in T(C)$, are weak-* continuous and affine.
	For a finite convex combination of extreme traces, the triangle inequality gives $|\tau(u_N^d)|\le\varepsilon/2$.
	The nonnegative affine function given by the second evaluation is at most $\varepsilon^2/4$
	on the same convex combinations.
	The Krein--Milman theorem and continuity extend both inequalities to every trace in $T(C)$.
	Taking square roots in the second one gives \eqref{eq:unitarytarget}, with the strict bound
	$\varepsilon$ in both cases.
	Setting $u=u_N$ and $a=ta_N$ for one such $N$ finishes the proof.
\end{proof}

\subsection{Local equivariant property Gamma}\label{sec:gamma}

We are now ready to prove local equivariant property Gamma.

\begin{thm}\label{thm:gamma}
Let $B$ be a separable, unital, $\mathcal{Z}$-stable C*-algebra with $T(B)\ne\emptyset$.
Let $\beta\colon G\to\Aut(B)$ be an action of a countable discrete group, and let
$E\subseteq Z(B)$ be a unital $G$-invariant C*-subalgebra on which $\beta$ is amenable.
Then $\beta$ has local equivariant property Gamma.
\end{thm}
	\begin{proof}
	Let $S\subseteq B^\omega$ be uniform $2$-norm separable, and fix an integer $m\ge2$.
	Choose a sequence $(s_j)_{j\in\N}$ of positive contractions whose linear span contains a
	uniform $2$-norm dense subset of $S$ and which includes the unit.
	Let $L$ be the separable $G$-invariant C*-subalgebra of $B^\omega$ generated by
	$\Delta_B(B)$ and the elements $s_j$, for $j\in\N$.
	Use \autoref{prop:rows} to find commuting $G$-invariant rows $C_r$ in $B^\omega\cap L'$ for all $r\in\N$.
	Define the elements $x_{r,g}\in C_r$ as in \autoref{rem:row-elements}, and set
	$C_0=C^*(\bigcup_{r\in\N}C_r)$ and $C=C^*(E,C_0)$.
	The algebra $C$ commutes with every $s_j$, since the image of $E$ is central in $B^\omega$
	and the rows commute with $L$.
	For every $\tau\in T_\omega(B)$, the restriction $\tau|_C$ is a tracial state on $C$, so $T(C)\ne\emptyset$.
	Moreover, the normalizations \eqref{eq:rowmoments} hold under every trace on $C_0$.
	Applying $\Delta_B$ to the weights in \eqref{eq:weights} and using its contractivity shows
	that the image of $E$ in $B^\omega$ has amenability weights.
	Fix a finite subset $K\subseteq G$, an integer $k\ge1$, and $\varepsilon>0$.
	Apply \autoref{prop:unitaries} to the restricted action on $C$, with $m=k$, and choose the
	resulting unitary $u\in C$.
	Thus $|\sigma(u^d)|<\varepsilon$ and $\|\beta_h^\omega(u)-u\|_{2,\sigma}<\varepsilon$ for
	every $\sigma\in T(C)$, every integer $d$ with $1\le|d|\le k$, and every $h\in K$.

	Let $\tau\in T_\omega(B)$ and let $j\in\N$ satisfy $\tau(s_j)>0$.
	Set $\sigma(c)=\tau(s_jc)/\tau(s_j)$ for all $c\in C$.
	Since $s_j$ commutes with $C$, the functional $\sigma$ is a tracial state on $C$.
	It follows that the chosen unitary $u$ satisfies
	\begin{equation}\label{eq:relativeapprox}
 |\tau(s_ju^d)|\le\varepsilon\tau(s_j)\le\varepsilon
	\end{equation}\nopagebreak
	for all $j\in\N$, all $d\in\Z$ with $1\le|d|\le k$, and all $\tau\in T_\omega(B)$; when
	$\tau(s_j)=0$, the same inequality follows from positivity.
	Restricting each limit trace to $C$ also gives
	\[
 \sup_{\tau\in T_\omega(B)}\|\beta_h^\omega(u)-u\|_{2,\tau}\le\varepsilon
	\]\nopagebreak
	for all $h\in K$.

	For every $j\in\N$, fix a sequence $(s_{j,n})_{n\in\N}$ of positive contractions in $B$ representing $s_j$.
	Choose an increasing sequence $(K_k)_{k\in\N}$ of finite subsets of $G$ with union equal to $G$.
	For each $k\in\N$, apply \autoref{prop:unitaries} with $K=K_k$, $m=k$, and $\varepsilon=\frac{1}{2k}$.
	Write the resulting unitary as $u^{(k)}=e^{ia^{(k)}}$.
	Choose a bounded selfadjoint representative $(a_n^{(k)})_{n\in\N}$ of $a^{(k)}$, and set
	$u_n^{(k)}=e^{ia_n^{(k)}}$.
	Then $(u_n^{(k)})_{n\in\N}$ is a representative consisting of unitaries in $B$.

	At coordinate $n$, take the set $X_n$ in \autoref{lem:diagonal} to be the unitary group of $B$.
	The functions to which we apply the lemma are given by the following expressions, evaluated at $w\in X_n$:
	\[
 \uu{[w,s_{j,n}]},\quad
 \uu{\beta_h(w)-w},
 \quad\text{and}\quad
 \sup_{\nu\in T(B)}|\nu(s_{j,n}w^d)|
	\]\nopagebreak
	for $j\in\N$, $h\in G$, and $d\in\Z\setminus\{0\}$.
	Since $u^{(k)}$ belongs to $C$, it commutes with every $s_j$.
	Thus \eqref{eq:ultranorm} gives
	\[
 \lim_{n\to\omega}\uu{[u_n^{(k)},s_{j,n}]}=0
 \quad\text{and}\quad
 \lim_{n\to\omega}\uu{\beta_h(u_n^{(k)})-u_n^{(k)}}\le \frac{1}{2k}
	\]\nopagebreak
	for all $j\in\N$ and all $h\in K_k$.
	For the trace expressions, we have
	\[
 \lim_{n\to\omega}\sup_{\nu\in T(B)}|\nu\bigl(s_{j,n}(u_n^{(k)})^d\bigr)|
 \stackrel{\eqref{eq:trace-sup}}{=}\sup_{\tau\in T_\omega(B)}|\tau\bigl(s_j(u^{(k)})^d\bigr)|
 \stackrel{\eqref{eq:relativeapprox}}{\le}\frac1{2k}.
	\]\nopagebreak
	for every $j\in\N$ and every $d\in\Z$ with $1\le|d|\le k$.
	Given finitely many of these functions and a positive tolerance, choose $k$ so large that
	$K_k$ contains the group elements involved, $k$ bounds the absolute values of the powers
	involved, and $\frac{1}{2k}$ is smaller than the tolerance.
	The sequence $(u_n^{(k)})_{n\in\N}$ then satisfies all the required inequalities.
	By \autoref{lem:diagonal}, there is a sequence $(u_n)_{n\in\N}$ of unitaries in $B$ on
	which every one of these functions has ultralimit zero.
	Their class $u=[(u_n)_{n\in\N}]$ is a unitary in
	$(B^\omega)^{\beta^\omega}\cap\{s_j\colon j\in\N\}'$ and satisfies
	\begin{equation}\label{eq:Haarmoments}
 \tau(s_ju^d)=0
	\end{equation}\nopagebreak
	for all $j\in\N$, all $d\in\Z\setminus\{0\}$, and all $\tau\in T_\omega(B)$.
	Denote by $\mu$ Haar probability measure on $\T$.
	Then \eqref{eq:Haarmoments} gives
	\begin{equation}\label{eq:Haarcontinuous}
 \tau(s_jf(u))=\tau(s_j)\int_{\T}f\,d\mu
	\end{equation}\nopagebreak
	for every trigonometric polynomial $f$, every $j\in\N$, and every $\tau\in T_\omega(B)$.
	By density and boundedness, \eqref{eq:Haarcontinuous} also holds for every $f\in C(\T)$.

	We finally construct the projections required by \autoref{def:gamma}, using the
	approximation by continuous functions explained in \autoref{rem:spectral-projections}.
	Partition $\T$ into $m$ half open arcs $I_1,\ldots,I_m$ of Haar measure $\frac{1}{m}$.
	For each $r=1,\ldots,m$, choose a sequence $(f_{r,k})_{k\in\N}$ of continuous functions
	from $\T$ to $[0,1]$ such that
	\[
 \lim_{k\to\infty}\|f_{r,k}-\indicator_{I_r}\|_{L^2(\mu)}=0
	\]\nopagebreak
	for all $r=1,\ldots,m$.
	Since the unit was included among the $s_j$, the identity in \eqref{eq:Haarcontinuous} applied to $|f|^2$ gives
	\begin{equation}\label{eq:Haar2norm}
 \uu{f(u)}^2=\int_{\T}|f|^2\,d\mu
	\end{equation}\nopagebreak
	for every $f\in C(\T)$.
	Thus $(f_{r,k}(u))_{k\in\N}$ is Cauchy in uniform $2$-norm, and completeness of bounded
	balls gives a positive contraction $p_r\in B^\omega$ as its limit.

	The functions $f_{r,k}$ are uniformly bounded, so their $L^2$ convergence also gives
	\[
 \lim_{k\to\infty}\|f_{r,k}^2-f_{r,k}\|_{L^2}=0,
 \quad
 \lim_{k\to\infty}\|f_{r,k}f_{s,k}\|_{L^2}=0,
 \quad\text{and}\quad
 \lim_{k\to\infty}\Bigl\|\sum_{r=1}^m f_{r,k}-1\Bigr\|_{L^2}=0
	\]\nopagebreak
	for all distinct $r,s\in\{1,\ldots,m\}$.
	Using \eqref{eq:Haar2norm} and continuity of multiplication on bounded balls, these
	identities imply $p_r^2=p_r$, $p_rp_s=0$ for $r\ne s$, as well as $\sum_{r=1}^m p_r=1$.
	Since the action is isometric and commutators are continuous on bounded balls by
	\autoref{rem:ultranorm}, invariance and commutation with the $s_j$ pass from $f_{r,k}(u)$
	to its uniform $2$-norm limit $p_r$.

	By the Cauchy--Schwarz inequality, we have
	\[
 \bigl|\tau\bigl(s_j\bigl(p_r-f_{r,k}(u)\bigr)\bigr)\bigr|\le\uu{p_r-f_{r,k}(u)}
	\]\nopagebreak
	for all $j,k\in\N$, all $r=1,\ldots,m$, and all $\tau\in T_\omega(B)$.
	We can therefore pass to the limit in \eqref{eq:Haarcontinuous}, obtaining
	\[
 \tau(s_jp_r)=\frac{1}{m}\tau(s_j)
	\]\nopagebreak
	for all $j\in\N$, all $r=1,\ldots,m$, and all $\tau\in T_\omega(B)$.
	Linearity gives the identity and commutation for every element in the linear span of the $s_j$.
	For $s\in S$ and $s'$ in this span, we obtain
	\[
 \uu{[p_r,s]}\le2\uu{s-s'}
 \quad\text{and}\quad
 \Bigl|\tau(sp_r)-\frac{1}{m}\tau(s)\Bigr|\le\Bigl(1+\frac{1}{m}\Bigr)\uu{s-s'}
	\]\nopagebreak
	for all $r=1,\ldots,m$ and all $\tau\in T_\omega(B)$.
	Since this span is dense in $S$ for the uniform $2$-norm, the right hand sides can be made arbitrarily small.
	Thus $p_r$ commutes with every $s\in S$ and satisfies $\tau(sp_r)=\frac{1}{m}\tau(s)$ for
	every $s\in S$ and $\tau\in T_\omega(B)$.
	The projections $p_1,\ldots,p_m$ are invariant, pairwise orthogonal, and add to the unit,
	so they satisfy \autoref{def:gamma}.
\end{proof}

\section{Averaging, matrix algebras, and absorption}\label{sec:averaging}

In this section, we prove the main results.
We first establish an averaging lemma for central amenability weights, which allows us to apply
the tracial local-to-global principle of \cite{SWGamma}.
To obtain the local solutions needed for this principle, we use amenability and McDuff
absorption to solve the matrix unit relations in the orbit von Neumann algebras.
We then obtain invariant central matrix units for the auxiliary action, and use a final
diagonal argument to obtain the corresponding matrix units for the original action.

\subsection{Central averaging and the local-to-global principle}

In addition to local equivariant property Gamma, the proof of the local-to-global theorem in
\cite{SWGamma} uses an averaging lemma for amenable groups, for which we now prove an analogue
using central amenability weights.

We first recall the weak complemented partitions of unity that we will use before averaging over $G$.

\begin{prop}[{\cite[Proposition 3.1]{SWGamma}}]\label{prop:weakcpou}
Let $B$ be a unital C*-algebra with $T(B)\ne\emptyset$, let $S\subseteq B^\omega$ be uniform
$2$-norm separable, and let $a_1,\ldots,a_k\in B_+$.
Let $\delta>0$ satisfy
\[
 \sup_{\tau\in T(B)}\min_{1\le j\le k}\tau(a_j)<\delta.
\]
Then there are positive contractions $e_1,\ldots,e_k\in B^\omega\cap S'$ such that
\begin{equation}\label{eq:weakcpou}
 \tau\Bigl(\sum_{j=1}^k e_j\Bigr)=1
 \quad\text{and}\quad
 \tau(a_je_j)\le\delta\tau(e_j)
\end{equation}\nopagebreak
for all $\tau\in T_\omega(B)$ and all $j=1,\ldots,k$.
\end{prop}

To use these trace inequalities after averaging, we will change a limit trace by a central positive weight.
The following lemma shows that this operation still gives a limit trace.

\begin{lma}\label{lem:weightedtrace}
Let $B$ be a unital C*-algebra with $T(B)\ne\emptyset$, let $\tau\in T_\omega(B)$, and let
$z\in Z(B)_+$ satisfy $\tau(z)>0$.
Then the functional $\sigma\colon B^\omega\to\C$ given by
$\sigma(c)=\tau(zc)/\tau(z)$ for all $c\in B^\omega$ is a limit trace.
If $\beta\colon G\to\Aut(B)$ is an action and $g\in G$, the same holds for
the functional $\sigma_g\colon B^\omega\to\C$ given by
$\sigma_g(c)=\tau(z\beta_g^\omega(c))/\tau(z)$ for all $c\in B^\omega$.
\end{lma}
	\begin{proof}
	Choose a representing sequence $(\tau_n)_{n\in\N}$ in $T(B)$.
	Set
	\[
 N_\tau=\bigl\{n\in\N\colon\tau_n(z)>\tfrac{1}{2}\tau(z)\bigr\},
	\]
	and note that $N_\tau$ belongs to $\omega$.
	For $n\in N_\tau$, set $\nu_n(b)=\tau_n(zb)/\tau_n(z)$ for all $b\in B$.
	For $n\notin N_\tau$, set $\nu_n=\tau_n$.
	Centrality and positivity of $z$ imply that each $\nu_n$ is a tracial state on $B$.
	For every bounded sequence $(b_n)_{n\in\N}$, we have
	\[
 \lim_{n\to\omega}\nu_n(b_n)
 =\frac{\lim_{n\to\omega}\tau_n(zb_n)}{\lim_{n\to\omega}\tau_n(z)}.
	\]
	Thus the sequence $(\nu_n)_{n\in\N}$ induces the limit trace $\sigma$.
	Likewise, $(\nu_n\circ\beta_g)_{n\in\N}$ induces $\sigma_g$, proving the second assertion.
\end{proof}

We also recall the following standard fact.

\begin{rem}\label{rem:trace-one-norm}
Let $D$ be a unital C*-algebra and let $\tau\in T(D)$.
For $x\in D$, set $\|x\|_{1,\tau}=\tau(|x|)$.
This defines a seminorm on $D$; in particular,
\[
 \tau(|x+y|)\le\tau(|x|)+\tau(|y|)
\]
for all $x,y\in D$.
We use this notation also for $D=B^\omega$ and $\tau\in T_\omega(B)$.
\end{rem}

\begin{lma}\label{lem:averaging}
Let $B$ be a unital C*-algebra with $T(B)\ne\emptyset$, and let $\beta\colon G\to\Aut(B)$ be an
action of a countable discrete group.
Suppose that $E\subseteq Z(B)$ is a unital $G$-invariant C*-subalgebra on which the action is amenable.
For every $\varepsilon>0$ and finite $K\subseteq G$, there is a finite nonempty set
$L\subseteq G$ with the following property.
For every uniform $2$-norm separable subset $S\subseteq B^\omega$, every
$a_1,\ldots,a_k\in B_+$, and every $\delta>0$ satisfying
\begin{equation}\label{eq:threshold}
 \sup_{\tau\in T(B)}\min_{1\le j\le k}\tau(a_j)<\frac{\delta}{|L|},
\end{equation}
there are positive contractions $e_1,\ldots,e_k\in B^\omega\cap S'$ such that
\begin{align}
 \tau\Bigl(\sum_{j=1}^ke_j\Bigr)&=1,\label{eq:averagesum}\\
 \tau(a_je_j)&\le\delta\tau(e_j),\label{eq:averagetrace}\\
 \sum_{j=1}^k\|\beta_h^\omega(e_j)-e_j\|_{1,\tau}&<\varepsilon\label{eq:averageinv}
\end{align}\nopagebreak
for all $\tau\in T_\omega(B)$, all $j=1,\ldots,k$, and all $h\in K$.
\end{lma}
	\begin{proof}
	Choose amenability weights $(d_g)_{g\in G}$ in $E$ satisfying \eqref{eq:weights} for $K$
	and $\varepsilon$, and set $L=\{g\in G\colon d_g\ne0\}$.
	Let $S$, $a_1,\ldots,a_k$, and $\delta$ be as in the statement.

	For data satisfying \eqref{eq:threshold}, set
	$a_j'=\frac{1}{|L|}\sum_{g\in L}\beta_{g^{-1}}(a_j)$ for $j=1,\ldots,k$.
	For each $\tau\in T(B)$, the functional $\frac{1}{|L|}\sum_{g\in L}\tau\circ\beta_{g^{-1}}$
	is a tracial state, and therefore
	\[
 \sup_{\tau\in T(B)}\min_{1\le j\le k}\tau(a_j')
 \le\sup_{\nu\in T(B)}\min_{1\le j\le k}\nu(a_j)
 \stackrel{\eqref{eq:threshold}}{<}\frac{\delta}{|L|}.
	\]
	Set $T=\bigcup_{g\in G}\beta_g^\omega(S)$, which is uniform $2$-norm separable since $G$ is countable.
	Apply \autoref{prop:weakcpou} to $T$, the positive elements $a_j'$, and
	$\frac{\delta}{|L|}$ in place of $\delta$, to obtain positive contractions
	$e_j'\in B^\omega\cap T'$ such that
	\begin{equation}\label{eq:prime}
 \tau\Bigl(\sum_{j=1}^k e_j'\Bigr)=1
 \quad\text{and}\quad
 \tau(a_j'e_j')\le\frac{\delta}{|L|}\tau(e_j')
	\end{equation}\nopagebreak
	for all $\tau\in T_\omega(B)$ and all $j=1,\ldots,k$.
	For every $g\in L$, positivity gives $\beta_{g^{-1}}(a_j)\le |L|a_j'$.
	Applying the positive functional defined by
	$f(a)=\tau\bigl((e_j')^{\tfrac{1}{2}}a(e_j')^{\tfrac{1}{2}}\bigr)$ and using traciality
	gives
	\begin{equation}\label{eq:individual}
 \tau(\beta_{g^{-1}}(a_j)e_j')
 \le |L|\tau(a_j'e_j')
 \stackrel{\eqref{eq:prime}}{\le}\delta\tau(e_j')
	\end{equation}\nopagebreak
	for all $g\in L$, all $\tau\in T_\omega(B)$, and all $j=1,\ldots,k$.

	For $j=1,\ldots,k$, set $e_j=\sum_{g\in L}d_g\beta_g^\omega(e_j')$.
	Since the $d_g$ are central positive elements with sum one, we have $0\le e_j\le1$.
	For $s\in S$, the element $e_j'$ commutes with $\beta_{g^{-1}}^\omega(s)$, so
	$\beta_g^\omega(e_j')$ commutes with $s$.
	Hence $e_j\in S'$.

	Fix $\tau\in T_\omega(B)$ and let $g\in L$.
	Suppose first that $\tau(d_g)>0$.
	By \autoref{lem:weightedtrace}, the functional
	$\sigma_g(c)=\tau(d_g\beta_g^\omega(c))/\tau(d_g)$, for $c\in B^\omega$, is a limit trace, and thus
	\begin{align*}
 \tau(a_jd_g\beta_g^\omega(e_j'))
 &=\tau(d_g)\sigma_g(\beta_{g^{-1}}(a_j)e_j')\\
 &\stackrel{\eqref{eq:individual}}{\le}
 \delta\tau(d_g)\sigma_g(e_j')
 =\delta\tau(d_g\beta_g^\omega(e_j')).
	\end{align*}
	If $\tau(d_g)=0$, both expressions vanish, since the positive functional with value
	$\tau(d_gx)$ at $x$ has norm $\tau(d_g)$.
	Adding over $g$ proves \eqref{eq:averagetrace}.
	Applying the same weighted traces to the first identity in \eqref{eq:prime} gives
	\[
 \tau\Bigl(d_g\beta_g^\omega\Bigl(\sum_{j=1}^k e_j'\Bigr)\Bigr)=\tau(d_g).
	\]
	Adding over $g\in L$ then proves \eqref{eq:averagesum}.

	For $h\in K$ and $r\in G$, set $\Delta_{h,r}=\beta_h(d_{h^{-1}r})-d_r$, and recall that $d_g=0$ for $g\notin L$.
	Applying $\beta_h^\omega$ to the definition of $e_j$ and setting $r=hg$, we get
	\begin{align*}
 \beta_h^\omega(e_j)-e_j
 &=
   \sum_{g\in G}\beta_h(d_g)\beta_{hg}^\omega(e_j')
     -\sum_{r\in G}d_r\beta_r^\omega(e_j')\\
 &=\sum_{r\in G}\bigl(\beta_h(d_{h^{-1}r})-d_r\bigr)\beta_r^\omega(e_j')
 =\sum_{r\in G}\Delta_{h,r}\beta_r^\omega(e_j')
	\end{align*}\nopagebreak
	for all $h\in K$ and all $j=1,\ldots,k$.
	Since each $\Delta_{h,r}$ is central and selfadjoint, we have
	$|\Delta_{h,r}\beta_r^\omega(e_j')|=|\Delta_{h,r}|\beta_r^\omega(e_j')$.
	If $\tau(|\Delta_{h,r}|)>0$, define $\nu_{h,r}\colon B^\omega\to\C$ by
	\[
 \nu_{h,r}(c)=\frac{\tau(|\Delta_{h,r}|\beta_r^\omega(c))}{\tau(|\Delta_{h,r}|)}
	\]\nopagebreak
	for all $c\in B^\omega$.
	By \autoref{lem:weightedtrace}, $\nu_{h,r}$ is a limit trace, so the first identity in
	\eqref{eq:prime} gives $\nu_{h,r}(\sum_{j=1}^k e_j')=1$.
	Multiplying by $\tau(|\Delta_{h,r}|)$, we obtain
	\begin{equation}\label{eq:weightednormalization}
 \sum_{j=1}^k\tau(|\Delta_{h,r}|\beta_r^\omega(e_j'))=\tau(|\Delta_{h,r}|).
	\end{equation}
	The same equality holds when $\tau(|\Delta_{h,r}|)=0$, since each term on the left then vanishes.
	Using \autoref{rem:trace-one-norm} at the first step, we get
	\begin{align*}
 \sum_{j=1}^k\|\beta_h^\omega(e_j)-e_j\|_{1,\tau}
 &\le\sum_{r\in G}\sum_{j=1}^k\tau(|\Delta_{h,r}|\beta_r^\omega(e_j'))\\
 &\stackrel{\eqref{eq:weightednormalization}}{=}\sum_{r\in G}\tau(|\Delta_{h,r}|)
 \le\Bigl\|\sum_{r\in G}|\Delta_{h,r}|\Bigr\|
 \stackrel{\eqref{eq:weights}}{<}\varepsilon.
	\end{align*}
	This proves \eqref{eq:averageinv}.
\end{proof}

We now state the local-to-global principle in the form needed below.
A \emph{$G$-polynomial} is a finite noncommutative polynomial in variables, their adjoints, and
their formal $G$-translates, with coefficients in the indicated algebra.
Its evaluation uses the given action on the variables and coefficients.

\begin{thm}\label{thm:localglobal}
Let $B$ be a separable, unital, $\mathcal{Z}$-stable C*-algebra with $T(B)\ne\emptyset$.
Let $\beta\colon G\to\Aut(B)$ be an action of a countable discrete group, and let
$E\subseteq Z(B)$ be a unital $G$-invariant C*-subalgebra on which $\beta$ is amenable.
Let $(P_j)_{j\in\N}$ be a countable family of $G$-polynomials with coefficients in $B^\omega$
and countably many contraction variables.
Suppose that, for every $\tau\in\overline{T_\omega(B)}^{\,w^*}$, the equations $P_j=0$, for
$j\in\N$, have a simultaneous solution by elements of norm at most one in
$(N_\tau^{\beta^\omega})^\omega$, with coefficients mapped by the orbit representation.
Then these equations have a simultaneous solution by elements of norm at most one in $B^\omega$.
\end{thm}
	\begin{proof}
	By \autoref{thm:gamma}, the action $\beta$ has local equivariant property Gamma, and
	\autoref{prop:weakcpou} gives weak complemented partitions of unity for $B$.
	Note that \autoref{lem:averaging} gives the conclusion of \cite[Lemma 3.4]{SWGamma}, with
	the support set chosen before the positive elements and with $\delta/|L|$ in
	\eqref{eq:threshold}.
	By \cite[Remark 4.8]{SWGamma}, this conclusion and local equivariant property Gamma imply
	its Theorem 4.2 for actions of arbitrary countable groups.
	In its unital form, the conclusion of that theorem says that it suffices to find, for every
	$\tau\in\overline{T_\omega(B)}^{\,w^*}$, every $m\in\N$, and every $\varepsilon>0$,
	contractions in $B^\omega$ for which the first $m$ polynomial expressions have $2$-seminorm
	less than $\varepsilon$ with respect to $\tau$.
	We verify this condition from the hypothesis in the statement.

	Fix $\tau$, $m$, and $\varepsilon$ as above, and write $N=N_\tau^{\beta^\omega}$ and
	$\pi=\pi_\tau^{\beta^\omega}$.
	Use the faithful normal trace $\widehat\tau$ from \autoref{rem:orbittrace}, and choose
	$\varepsilon_0>0$ such that $2\varepsilon_0<\sqrt{c_e}\varepsilon$.
	For each $j\in\N$, write $P_j^N$ for the polynomial obtained by representing the
	coefficients of $P_j$ using $\pi$ and extending the action normally to $N$.
	Only finitely many variables and group elements occur in $P_1,\ldots,P_m$.
	Choose representatives of norm at most one in $N$ for those entries of the solution in $N^\omega$.
	For $n\in\N$, let $\mathbf y_n$ be the tuple of these representatives at coordinate $n$.
	Then $\|P_j^N(\mathbf y_n)\|_{2,\widehat\tau}$ tends to zero along $\omega$ for each $j=1,\ldots,m$.
	Intersecting the corresponding finitely many $\omega$-large sets gives one coordinate, with
	a tuple $\mathbf y$ of elements of norm at most one in $N$, such that
	$\|P_j^N(\mathbf y)\|_{2,\widehat\tau}<\varepsilon_0$ for all $j=1,\ldots,m$.

	By Kaplansky density, each entry of $\mathbf y$ can be approximated in the strong-*
	topology by elements of norm at most one in $\pi(B^\omega)$.
	Recall that each such element has a lift of norm at most one in $B^\omega$, since $\pi$ is a homomorphism.
	On bounded sets, taking adjoints, multiplying finitely many factors, and applying each of
	the finitely many normal automorphisms are continuous for the strong-* topology.
	Choose the approximations close enough that their lifts form a tuple $\mathbf z$ of
	elements of norm at most one satisfying
	\[
 \bigl\|P_j^N(\pi(\mathbf z))-P_j^N(\mathbf y)\bigr\|_{2,\widehat\tau}<\varepsilon_0
	\]\nopagebreak
	for all $j=1,\ldots,m$.
	Equivariance and multiplicativity of $\pi$ give $P_j^N(\pi(\mathbf z))=\pi(P_j(\mathbf z))$.
	Therefore
	\[
 \|P_j(\mathbf z)\|_{2,\tau}
 \stackrel{\eqref{eq:orbittrace}}{\le}
 \frac{1}{\sqrt{c_e}}\bigl\|P_j^N(\pi(\mathbf z))\bigr\|_{2,\widehat\tau}
 <\frac{2\varepsilon_0}{\sqrt{c_e}}<\varepsilon
	\]\nopagebreak
	for every $j=1,\ldots,m$.

	This proves the local approximation condition, and the local-to-global theorem gives
	contractions in $B^\omega$ satisfying all the relations.
\end{proof}

\subsection{Invariant central matrix algebras}\label{sec:matrix}

The local solutions required by \autoref{thm:localglobal} come from the McDuff absorption
theorem recalled in \autoref{sec:prelim}.
We will apply it to the orbit von Neumann algebra of each trace on $B$, and then regard the
resulting matrix units as local solutions for the uniform tracial ultrapower.

We will use the following standard inequality for central coefficients.
\begin{rem}\label{rem:central-sum-bound}
Let $F$ be a finite set, let $M$ be a von Neumann algebra, let $(a_g)_{g\in F}$ be selfadjoint
elements of $Z(M)$, and let $(y_g)_{g\in F}$ be contractions in $M$.
Then
\[
 \Bigl\|\sum_{g\in F}a_gy_g\Bigr\|\le\Bigl\|\sum_{g\in F}|a_g|\Bigr\|.
\]
\end{rem}

\begin{prop}\label{prop:matricesB}
Let $B$ be a separable, unital, $\mathcal{Z}$-stable C*-algebra with $T(B)\ne\emptyset$.
Let $\beta\colon G\to\Aut(B)$ be an action of a countable discrete group, and let
$E\subseteq Z(B)$ be a unital $G$-invariant C*-subalgebra on which $\beta$ is amenable.
For every $m\in\N$, there is a unital homomorphism
\[
 M_m\to(B^\omega\cap B')^{\beta^\omega}.
\]
\end{prop}
	\begin{proof}
	The assertion is immediate for $m=1$.
	For $m\ge2$, we first construct local invariant matrix units, and then check that they
	solve the equations required by \autoref{thm:localglobal}.
	Fix an integer $m\ge2$, and write $(e_{st})_{1\le s,t\le m}$ for the standard matrix units of $M_m$.
	Fix $\sigma\in T(B)$, and set $M_\sigma=N_\sigma^\beta$, using the orbit von Neumann
	algebra from \autoref{def:orbit}.
	Its faithful normal tracial state restricts on $B$ to
	$\tau_\sigma=\sum_{g\in G}c_g\sigma\circ\beta_{g^{-1}}$.
	Since the represented copy of $B$ is weakly dense, $M_\sigma$ is canonically isomorphic to
	$\pi_{\tau_\sigma}(B)''$, and hence is McDuff and has separable predual by
	\autoref{lem:gns}.
\vspace{.2cm}

	\textbf{Claim 1:} \emph{The induced action on $M_\sigma$ is amenable.} Write $\gamma$ for this action and continue to denote the represented weights by $d_g$.
	Since $E$ commutes with $B$, its image commutes with the weak closure $M_\sigma$ and lies in $Z(M_\sigma)$.
	The represented weights satisfy \eqref{eq:weights}, by contractivity of the representation.
	Define $\Phi_d\colon\ell^\infty(G)\bar\otimes M_\sigma\to M_\sigma$ by
	\[
 \Phi_d(x)=\sum_{g\in G}d_gx_g
	\]\nopagebreak
	for every bounded family $x=(x_g)_{g\in G}$ in $M_\sigma$.
	The map $\Phi_d$ is unital since the weights add to one, and completely positive since
	every $d_g$ is positive and central.
	Centrality also implies that $\Phi_d$ is an $M_\sigma$-bimodule map.
	Moreover, $\Phi_d(x)=a$ whenever $x_g=a$ for all $g\in G$ and some $a\in M_\sigma$.

	Recall that the diagonal action on the domain satisfies
	$(h\cdot x)_g=\gamma_h(x_{h^{-1}g})$ for all $g,h\in G$ and all
	$x\in\ell^\infty(G)\bar\otimes M_\sigma$.
	Reindexing the sum for $\gamma_h(\Phi_d(x))$ gives
	\[
 \Phi_d(h\cdot x)-\gamma_h\bigl(\Phi_d(x)\bigr)
 =\sum_{g\in G}(d_g-\gamma_h(d_{h^{-1}g}))\gamma_h(x_{h^{-1}g}).
	\]
	By \autoref{rem:central-sum-bound}, we have
	\begin{equation}\label{eq:expectation-defect}
 \bigl\|\Phi_d(h\cdot x)-\gamma_h\bigl(\Phi_d(x)\bigr)\bigr\|
 \le\|x\|\Bigl\|\sum_{g\in G}|d_g-\gamma_h(d_{h^{-1}g})|\Bigr\|.
	\end{equation}
	Choose a sequence of families of weights for which the bounds in \eqref{eq:weights} tend to
	zero on every finite subset of $G$.
	For each $x$ in the domain, the values $\Phi_d(x)$ lie in the ultraweakly compact ball of radius $\|x\|$.
	Compactness of the product of these balls gives a subnet on which the maps converge
	pointwise in the ultraweak topology.
	Denote the limit map by $\Phi$.
	It is a unital completely positive $M_\sigma$-bimodule map fixing the constant copy of $M_\sigma$.
	Normality of $\gamma_h$ and \eqref{eq:expectation-defect} imply that $\Phi$ is equivariant.
	It follows that $\Phi\colon\ell^\infty(G)\bar\otimes M_\sigma\to M_\sigma$ is an
	equivariant conditional expectation, which proves the claim.
\vspace{.2cm}

	All the hypotheses of \autoref{thm:vn-absorption} are now verified, so the theorem gives a
	unital homomorphism $M_m\to(M_\sigma^\omega\cap M_\sigma')^{\gamma^\omega}$.
	Now fix $\tau\in\overline{T_\omega(B)}^{\,w^*}$ and set $\sigma=\tau|_B$.
	Set $N_\tau=N_\tau^{\beta^\omega}$, using the orbit von Neumann algebra of $B^\omega$.
\vspace{.2cm}

\textbf{Claim 2:} \emph{The homomorphism $\Delta_B$ induces a normal equivariant embedding $M_\sigma\to N_\tau$, and an equivariant embedding $M_\sigma^\omega\to N_\tau^\omega$.} To prove the claim, use the same strictly positive weights $c_g$ in the faithful normal traces on both algebras.
	The restriction of the trace of $N_\tau$ to $B$ is $\tau_\sigma$, the trace defining $M_\sigma$.
	Let $N_0$ be the von Neumann closure of the image of $B$ in $N_\tau$.
	For $b,c\in B$, the inner product of their images in $L^2(N_0,\widehat\tau)$ is
	\[
 \widehat\tau(\pi_\tau^{\beta^\omega}(c)^*\pi_\tau^{\beta^\omega}(b))
 =\tau_\sigma(c^*b).
	\]
	Since the images of $B$ are dense in this $L^2$-space by Kaplansky density, the map from
	the GNS space of $(B,\tau_\sigma)$ to $L^2(N_0,\widehat\tau)$ is unitary and intertwines
	the left representations of $B$.
	Taking von Neumann closures identifies $M_\sigma$ normally with $N_0$, and composing this
	identification with the inclusion of $N_0$ in $N_\tau$ gives the required normal embedding.
	The map is equivariant on the weakly dense copy of $B$, and hence on its closure.
	Moreover, it preserves the chosen faithful traces, so a bounded sequence is $2$-norm null
	along $\omega$ if and only if its image is, and we obtain an equivariant embedding of the
	tracial ultrapowers, proving the claim.
\vspace{.2cm}

	It follows that there is a unital homomorphism
	$M_m\to(N_\tau^\omega\cap B')^{\gamma^\omega}$, using $B$ for its represented image and
	$\gamma$ for the induced action on $N_\tau$.

	Choose a norm dense sequence $(b_j)_{j\in\N}$ in $B_1$.
	For variables $X_{st}$ of norm at most one, with $1\le s,t\le m$, consider the matrix unit relations
	\[
 X_{st}X_{uv}=\delta_{tu}X_{sv},\quad X_{st}^*=X_{ts},
 \quad\text{and}\quad \sum_{s=1}^m X_{ss}=1,
	\]
	together with $[X_{st},b_j]=0$ and $g\cdot X_{st}=X_{st}$ for all matrix indices in
	$\{1,\ldots,m\}$, all $j\in\N$, and all $g\in G$.
	The preceding construction solves this system in $N_\tau^\omega$ for every $\tau$ in the
	required weak-* closure. \autoref{thm:localglobal} gives solutions $y_{st}\in B^\omega$,
	for $1\le s,t\le m$.
	It is then easy to see that the unital homomorphism $\psi\colon M_m\to B^\omega$ determined
	by $\psi(e_{st})=y_{st}$ has range in $(B^\omega\cap B')^{\beta^\omega}$, as desired.
\end{proof}

\subsection{Transfer and absorption}\label{sec:main}

To complete the proof of the main results, it remains to pass from invariant central matrix
units for the auxiliary action to such matrix units for the original action.
Since the map between the two algebras takes values in $A^\omega$, this requires choosing
representatives for two successive ultrapowers.
We need only lift the images of individual matrix-unit candidates to contractions in $A$, and
not the homomorphism on the whole of $B$.
The resulting diagonal argument is standard, so we state the following lemma without proof.

\begin{lma}\label{lem:transfer}
Let $A$ and $B$ be separable unital C*-algebras with nonempty tracial state spaces.
Let $\alpha\colon G\to\Aut(A)$ and $\beta\colon G\to\Aut(B)$ be actions of a countable discrete group.
Let $\iota\colon A\to B$ and $\theta\colon B\to A^\omega$ be equivariant unital homomorphisms
satisfying $\theta\circ\iota=\Delta_A$.
Let $m\ge2$ be an integer.
If there is a unital homomorphism $M_m\to(B^\omega\cap\iota(A)')^{\beta^\omega}$, then there is
a unital homomorphism
\[
 M_m\to(A^\omega\cap A')^{\alpha^\omega}.
\]
\end{lma}

We are now ready to establish \autoref{mainThmB} from the introduction.

\begin{thm}\label{thm:tracial-matrices}
Let $A$ be a separable, unital, $\mathcal{Z}$-stable C*-algebra with $T(A)\ne\emptyset$.
Let $G$ be a countable discrete group, and let $\alpha\colon G\to\Aut(A)$ be tracially amenable.
For every integer $m\ge2$, there is a unital homomorphism
\[
 M_m\to (A^\omega\cap A')^{\alpha^\omega}.
\]
\end{thm}
	\begin{proof}
	Choose $B$, $\beta$, $E$, and $\theta$ as in \autoref{prop:auxiliary}.
	The inclusion $\iota\colon A\to B$ is given by $\iota(a)=a\otimes1$ for $a\in A$.
	For each $m\ge2$, \autoref{prop:matricesB} gives a unital homomorphism $M_m\to(B^\omega\cap B')^{\beta^\omega}$.
	Since $\theta\circ\iota=\Delta_A$, \autoref{lem:transfer} gives the required homomorphism
	into $(A^\omega\cap A')^{\alpha^\omega}$.
\end{proof}

Finally, we prove the main result of this paper, stated as \autoref{mainThmA} in the introduction.

\begin{thm}\label{thm:absorption}
Let $A$ be a separable, simple, unital, nuclear, $\mathcal{Z}$-stable C*-algebra.
Let $G$ be a countable discrete group, and let $\alpha\colon G\to\Aut(A)$ be an amenable action.
Then
\[
 (A,\alpha)\cc(A\otimes\mathcal{Z},\alpha\otimes\id_{\mathcal{Z}}).
\]
\end{thm}
	\begin{proof}
	Suppose first that $T(A)=\emptyset$.
	Since $A$ is simple, unital, nuclear, and $\mathcal{Z}$-stable, \cite[Corollary
	5.1]{Rordam} implies that $A$ is purely infinite, and hence a Kirchberg algebra.
	By \cite[Corollary 3.4]{SI}, we have $\alpha\cc\alpha\otimes\id_{\mathcal{O}_\infty}$.
	Tensoring this cocycle conjugacy with $\id_{\mathcal{Z}}$ and using
	$\mathcal{O}_\infty\otimes\mathcal{Z}\cong\mathcal{O}_\infty$, we obtain
	\[
 \alpha\otimes\id_{\mathcal{Z}}
 \cc\alpha\otimes\id_{\mathcal{O}_\infty}\otimes\id_{\mathcal{Z}}
 \cong\alpha\otimes\id_{\mathcal{O}_\infty}
 \cc\alpha.
	\]

	Suppose now that $T(A)\ne\emptyset$.
	Amenability of $\alpha$ implies tracial amenability, so \autoref{thm:tracial-matrices}
	gives invariant central matrix algebras of every size in $A^\omega$.
	Since $A$ is simple, nuclear, and $\mathcal{Z}$-stable, it has strict comparison by \cite[Corollary 4.6]{Rordam}.
	Hence $\alpha$ has equivariant property (SI) by \cite[Corollary 3.3]{SI}.
	The absorption criterion in \autoref{thm:si-absorption} now gives
	$\alpha\cc\alpha\otimes\id_{\mathcal{Z}}$, as desired.
\end{proof}

\begin{ex}\label{ex:suzuki}
Suzuki's construction in \cite[Theorem 1.3 and Remark 1.5]{Suzuki} shows that
\autoref{thm:absorption} applies to amenable actions of nonamenable groups on algebras with
traces, and not only to actions of amenable groups.
Starting with a free minimal action of a countably infinite locally finite group on a Cantor
space, he obtains a unital simple separable AF-algebra $A$ and amenable free group actions on
$A$.
The algebra $A$ is infinite dimensional, since it contains the continuous functions on the
Cantor space, and hence it is $\mathcal{Z}$-stable by \cite{JS}.
It is also nuclear and has traces, so \autoref{thm:absorption} implies that these actions are
equivariantly $\mathcal{Z}$-stable.
\end{ex}

\subsection{Crossed products and intermediate algebras}\label{sec:applications}
We conclude with the two applications stated in the introduction.
The first uses the fact that an invariant central copy of $\mathcal{Z}$ in the norm ultrapower
of $A$ remains central after passing to the crossed product.
This is also the argument underlying \cite[Proposition 5.1 and Corollary 4.5]{Sarkowicz}.

\begin{cor}\label{cor:intermediate}
Let $A$ be a separable, simple, unital, nuclear, $\mathcal{Z}$-stable C*-algebra, let $G$ be a
countable discrete group, and let $\alpha\colon G\to\Aut(A)$ be amenable.
Then every intermediate C*-algebra
\[
 A\subseteq C\subseteq A\rtimes_{\alpha,r}G
\]
is $\mathcal{Z}$-stable.
In particular, $A\rtimes_{\alpha,r}G$ is $\mathcal{Z}$-stable.
\end{cor}
	\begin{proof}
	Set $B=A\rtimes_{\alpha,r}G$, and denote by $u_g\in B$ the canonical unitary associated to $g\in G$.
	By \autoref{thm:absorption}, the action $\alpha$ is equivariantly $\mathcal{Z}$-stable, so
	the central sequence characterization of equivariant absorption gives a unital homomorphism
	\[
 \psi\colon\mathcal{Z}\to(A_\omega\cap A')^{\alpha_\omega};
	\]
	see also \cite[proof of Proposition 5.1]{Sarkowicz}.
	Moreover, the inclusion $A\subseteq B$ induces an isometric inclusion $A_\omega\subseteq B_\omega$.
	For $z\in\mathcal{Z}$ and $g\in G$, we have
	\[
 u_g\psi(z)u_g^*=\alpha_{\omega,g}(\psi(z))=\psi(z),
	\]
	so $\psi(z)$ commutes with each $u_g$ as well as with $A$.
	Since these elements generate $B$, the range of $\psi$ lies in $A_\omega\cap B'$.

	Now let $A\subseteq C\subseteq B$ be an intermediate C*-algebra.
	Then $A_\omega\subseteq C_\omega\subseteq B_\omega$, and the same map $\psi$ gives a unital
	homomorphism $\mathcal{Z}\to C_\omega\cap C'$.
	The algebra $C$ contains the unit of $B$ and is separable, so $C$ is $\mathcal{Z}$-stable by
	\cite[Theorem 2.2]{TW}.
\end{proof}

\begin{cor}\label{cor:pure}
Let $A$ be a separable, simple, unital, nuclear, $\mathcal{Z}$-stable C*-algebra, let $G$ be a
countable discrete nonamenable group, and let $\alpha\colon G\to\Aut(A)$ be amenable and outer.
Then $A\rtimes_{\alpha,r}G$ is a unital Kirchberg algebra.
\end{cor}
	\begin{proof}
	Set $B=A\rtimes_{\alpha,r}G$.
	The algebra $B$ is unital and separable, and it is simple because $A$ is simple and
	$\alpha$ is outer; see \cite{Kishimoto}.
	Amenability of $\alpha$ and nuclearity of $A$ imply that $B$ is nuclear by \cite[Theorem 4.5]{AD}.
	Moreover, \autoref{cor:intermediate} gives $\mathcal{Z}$-stability of $B$.

	We claim that $T(B)=\emptyset$.
	Indeed, if $\tau\in T(B)$, then $\sigma=\tau|_A$ is a tracial state on $A$ satisfying
	\[
 \sigma(\alpha_g(a))=\tau(u_gau_g^*)=\tau(a)=\sigma(a)
	\]\nopagebreak
	for all $a\in A$ and $g\in G$.
	Thus $\sigma$ is $G$-invariant, which contradicts nonamenability of $G$ by \cite[Lemma
	2.9]{GGKNV}, since amenability of $\alpha$ implies tracial amenability when
	$T(A)\ne\emptyset$.
	Therefore, the simple, unital, nuclear, $\mathcal{Z}$-stable algebra $B$ has no traces, so
	it is purely infinite by \cite[Corollary 5.1]{Rordam}.
\end{proof}

\appendix
\section{Estimates for Gaussian moments}\label{sec:probability-appendix}

The construction of approximately invariant unitaries using Gaussian variables in \cite{OAI}
provides the starting point for the proof of equivariant $\mathcal{Z}$-stability in this paper.
Its use requires several facts about moments and Gaussian quadratic forms which are less
familiar in the study of C*-algebras.
We give their proofs here, together with the probability background needed to follow them.
In particular, the two estimates in \autoref{prop:moment-estimates} are proved in full, so
their use in the main text does not require any result from \cite{OAI} as a black box.
The row counting argument and the use of Gaussian quadratic forms follow that preprint; the
purpose of this appendix is to present these tools in a form that can be read independently and
used in other C*-algebraic arguments.

The Gaussian variables arise from adding centered elements in the commuting row algebras
constructed in \autoref{sec:rows}, normalized by $\frac{1}{\sqrt N}$ for $N$ rows.
Independence under extreme traces makes terms involving a row only once vanish, while
normalization makes terms involving higher multiplicities negligible; the remaining terms are
precisely the pairings in the Gaussian moment formula.
Thus Gaussian variables describe the limiting tracial moments of these sums.
The Fourier transforms of the laws of their quadratic forms then give estimates for moments of
exponentials in the original algebra, producing the unitaries from which we obtain projections
of equal trace.

The argument has three parts.
First, sums of independent centered rows have approximately Gaussian moments, even though the
elements within each row need not commute.
Second, averaged covariance elements can be replaced by their scalar traces inside any fixed joint moment.
Finally, these moment calculations give estimates for exponentials of selfadjoint elements, by
comparison with Gaussian quadratic forms.
We begin with the correspondence between traces and scalar probability measures, which explains
why scalar probability is useful in this setting.
All the probability notation recalled in \autoref{sec:probability} is developed again below.

\subsection{Probability measures and traces}\label{app:probability}
Let $(\Omega,\mathcal F,\mathbb P)$ be a probability space.
A real random variable is a measurable function $X\colon\Omega\to\R$, and its \emph{law} is the
probability measure $\mu_X$ on $\R$ defined by $\mu_X(V)=\mathbb P(X^{-1}(V))$ for every Borel
set $V\subseteq\R$.
For every $\mu_X$-integrable function $f\colon\R\to\C$, integration with respect to the law gives
\[
 \mathbb E(f(X))=\int_\Omega f(X)\,d\mathbb P=\int_\R f\,d\mu_X.
\]
Thus all expectations involving only $X$ can be computed from $\mu_X$ without specifying the
original probability space.
We say that $X$ is \emph{centered} if $\mathbb E(X)=0$.
For a square integrable real random variable, set
$\operatorname{Var}(X)=\mathbb E\bigl((X-\mathbb E(X))^2\bigr)$ and
$\|X\|_{L^2}=\mathbb E(X^2)^{\tfrac{1}{2}}$.
For centered $X$, its variance is therefore $\mathbb E(X^2)=\|X\|_{L^2}^2$.

For a finite family $X_1,\ldots,X_m$, its joint law is defined in the same way on $\R^m$.
The variables are \emph{independent} if their joint law is the product of their individual laws.
If this is the case, Fubini's theorem gives
\begin{equation}\label{eq:app-independence}
 \mathbb E\Bigl(\prod_{j=1}^m f_j(X_j)\Bigr)
 =\prod_{j=1}^m\mathbb E(f_j(X_j))
\end{equation}
for bounded Borel functions $f_j\colon\R\to\C$, and also when each $f_j$ is integrable with
respect to the law of $X_j$.
It is this factorization that will enter the row argument.

For a probability measure $\mu$ on $\R$, its \emph{Fourier transform} is the function
$\widehat\mu\colon\R\to\C$ given by $\widehat\mu(t)=\int_\R e^{itx}\,d\mu(x)$ for all $t\in\R$.
We use this name for the function called the characteristic function in probability theory.
Fourier transforms determine probability measures, also on $\R^m$; see \cite[Theorem 26.2 and p.~382]{Billingsley}.

\begin{rem}\label{rem:app-spectral}
Recall that, for a state $\rho$ on a unital C*-algebra $C$ and a selfadjoint element $a\in C$,
the restriction of $\rho$ to $C^*(1,a)\cong C(\operatorname{sp}(a))$ is represented by a
probability measure $\mu_{a,\rho}$ on $\operatorname{sp}(a)\subseteq\R$ satisfying
\[
 \rho(f(a))=\int_{\operatorname{sp}(a)} f\,d\mu_{a,\rho}
\]
for all $f\in C(\operatorname{sp}(a))$.
We regard $\mu_{a,\rho}$ as a measure on $\R$ by assigning measure zero to the complement of
$\operatorname{sp}(a)$, and call it the \emph{spectral distribution} of $a$ under $\rho$.
In particular, $\rho(a^r)$ is its $r$th moment and $\rho(e^{ita})=\widehat{\mu_{a,\rho}}(t)$
for all integers $r\ge0$ and all $t\in\R$.
For a unitary $u\in C$, the same construction represents the restriction of $\rho$ to
$C^*(1,u)\cong C(\operatorname{sp}(u))$ by a probability measure on
$\operatorname{sp}(u)\subseteq\T$.

Several noncommuting selfadjoint elements need not have a joint probability measure.
We instead retain their ordered moments $\rho(a_{g_1}\cdots a_{g_d})$, for $g_1,\ldots,g_d\in J$.
Convergence to the moments of real random variables $(X_g)_{g\in J}$ means convergence of these
numbers to $\mathbb E(X_{g_1}\cdots X_{g_d})$ for every fixed ordered word.
\end{rem}

We will use both independence and commutation of the row algebras $C_r$.
More precisely, let $(C_r)_{r\in\N}$ be pairwise commuting unital C*-subalgebras of a unital C*-algebra $C$.
For $\tau\in T(C)$, the rows are \emph{independent with respect to $\tau$} if
\[
 \tau(a_1\cdots a_s)=\prod_{k=1}^s\tau(a_k)
\]
whenever $a_1,\ldots,a_s$ belong to distinct row algebras.
If the rows are independent, a word can be grouped according to its rows and then evaluated by
multiplying the traces of its row products, preserving the order within each row.
This is the C*-algebraic counterpart of \eqref{eq:app-independence}.

\begin{rem}\label{rem:app-circle}
Let $C$ be a unital C*-algebra, let $\rho$ be a state on $C$, let $a\in C$ be selfadjoint, let
$t\in\R$, and set $u=e^{ita}$.
The spectral measure of $u$ is obtained from $\mu_{a,\rho}$ by applying the continuous function
$h_t\colon\R\to\T$ given by $h_t(x)=e^{itx}$ for all $x\in\R$.
Explicitly,
\[
 \rho(f(u))=\int_\R f(e^{itx})\,d\mu_{a,\rho}(x)
\]
for every continuous function $f\colon\T\to\C$.
In particular, $\rho(u^k)=\widehat{\mu_{a,\rho}}(kt)$ for all $k\in\Z$, which explains why
bounds for Fourier transforms give bounds for unitary moments.
\end{rem}

\begin{df}\label{def:haar}
The \emph{Haar probability measure} on $\T$ is the measure $\mu$ satisfying
\[
 \int_\T f\,d\mu=\frac{1}{2\pi}\int_0^{2\pi}f(e^{it})\,dt
\]
for all $f\in C(\T)$.
\end{df}

Recall that all nonzero integer moments of Haar probability measure vanish, and this property
determines the measure because trigonometric polynomials are uniformly dense in $C(\T)$.
Thus a unitary whose nonzero moments vanish under every trace has the same spectral
distribution under all those traces.
In a uniform tracial ultrapower, \autoref{rem:spectral-projections} then gives projections of
equal trace by taking arcs of equal measure.
The main argument obtains exact vanishing by making successively larger finite sets of moments
small and then passing to the ultrapower.

\subsection{Gaussian vectors and their moments}\label{app:gaussians}

Recall that a standard real Gaussian variable is a real random variable $U$ with density
$\frac{1}{\sqrt{2\pi}}e^{-\tfrac{x^2}{2}}$.
Completing the square gives
\begin{equation}\label{eq:app-gaussian-mgf}
 \mathbb E(e^{tU})=\frac{1}{\sqrt{2\pi}}\int_\R e^{tx-\frac{x^2}{2}}\,dx=e^{\frac{t^2}{2}}
\end{equation}
for all $t\in\R$.
Odd moments vanish by symmetry, while differentiation at zero gives $\mathbb E(U^2)=1$ and $\mathbb E(U^4)=3$.

A family of real random variables is a \emph{centered real Gaussian vector} if its joint law is
that of a real linear image of finitely many independent standard real Gaussian variables.
For example, choose independent standard real Gaussian variables $\eta_1,\ldots,\eta_m$, a
nonempty finite set $J$, and real numbers $a_{g,k}$, for $g\in J$ and $k=1,\ldots,m$, and set
$\xi_g=\sum_{k=1}^m a_{g,k}\eta_k$ for all $g\in J$.
The family $(\xi_g)_{g\in J}$ is then a centered real Gaussian vector, with covariance entries
\[
 R_{g,h}=\mathbb E(\xi_g\xi_h)=\sum_{k=1}^m a_{g,k}a_{h,k}
\]
for all $g,h\in J$.
This is the Gram matrix of the coefficient vectors, so it is real symmetric and positive semidefinite.
Conversely, let $R$ be a real symmetric positive semidefinite matrix indexed by $J$, and let
$\eta=(\eta_g)_{g\in J}$ have independent standard real Gaussian coordinates.
Then $\xi=R^{\tfrac{1}{2}}\eta$ is a centered real Gaussian vector with covariance $R$.
No invertibility assumption on $R$ is needed.

For $t_g\in\R$, for $g\in J$, independence and \eqref{eq:app-gaussian-mgf} give
\begin{equation}\label{eq:app-vector-mgf}
 \mathbb E\Bigl(e^{\sum_{g\in J}t_g\xi_g}\Bigr)
 =\exp\Bigl(\frac{1}{2}\sum_{g,h\in J}t_gt_hR_{g,h}\Bigr).
\end{equation}
The Fourier transform of the joint law has the same formula with a minus sign in the exponent.
Since a probability measure is determined by its Fourier transform, the joint law of $\xi$ depends only on $R$.

It is often useful to express the construction using vectors in a Euclidean space.
Choose vectors $(v_g)_{g\in J}$ in a finite dimensional real Euclidean space $V$ with
$\langle v_g,v_h\rangle=R_{g,h}$ for all $g,h\in J$.
Such vectors are called \emph{Gram vectors} for $R$.
A standard Gaussian vector $\eta$ on $V$ has independent standard Gaussian coordinates in an
orthonormal basis; its distribution is unchanged by orthogonal changes of basis.
The family $(\langle\eta,v_g\rangle)_{g\in J}$ has covariance $R$.
If $R_{g,g}=1$, then $v_g$ is a unit vector, and if $R$ is singular, the Gram vectors are linearly dependent.

A \emph{pairing} of $\{1,\ldots,d\}$ is a partition of this set into subsets of size two.
We denote the set of pairings by $\mathcal P_2(d)$, and note that $\mathcal P_2(d)=\emptyset$ when $d$ is odd.

\begin{lma}\label{lem:app-wick}
Let $(\xi_g)_{g\in J}$ be a centered real Gaussian vector with covariance $R$.
For every even integer $d\ge2$, one has
\begin{equation}\label{eq:wick}
 \mathbb E(\xi_{g_1}\cdots\xi_{g_d})
 =\sum_{\mathcal P\in\mathcal P_2(d)}\prod_{\{a,b\}\in\mathcal P}R_{g_a,g_b}
\end{equation}
for all $g_1,\ldots,g_d\in J$.
All odd joint moments vanish.
\end{lma}
	\begin{proof}
	Fix $g_1,\ldots,g_d\in J$, with repetitions allowed, and assign a real variable $s_a$ to
	each position $a=1,\ldots,d$.
	Set $F(s_1,\ldots,s_d)=\mathbb E(e^{\sum_{a=1}^d s_a\xi_{g_a}})$.
	Differentiating once with respect to each variable multiplies the integrand by
	$\xi_{g_1}\cdots\xi_{g_d}$, and evaluating at zero removes the exponential.
	Thus
	\[
 \frac{\partial^d F}{\partial s_1\cdots\partial s_d}(0,\ldots,0)
 =\mathbb E(\xi_{g_1}\cdots\xi_{g_d}).
	\]
	These differentiations are justified by the exponential integrability of linear
	combinations of standard Gaussian variables.
	By \eqref{eq:app-vector-mgf}, we also have
	\[
 F(s_1,\ldots,s_d)=\exp\Bigl(\frac{1}{2}\sum_{a,b=1}^d s_as_bR_{g_a,g_b}\Bigr).
	\]
	The required derivative is the coefficient of $s_1\cdots s_d$ in this power series.
	Every term on the right has even total degree, so this coefficient is zero for odd $d$.
	If $d=2m$, only the $m$th power of the quadratic expression contributes.
	To obtain $s_1\cdots s_{2m}$, the $2m$ positions must be divided into $m$ pairs, and each
	pair $\{a,b\}$ contributes $R_{g_a,g_b}$.
	For a fixed pairing, there are $2^m$ choices of orientation of its pairs and $m!$ orders in which to choose them.
	These cancel the factor $1/(2^m m!)$ in the exponential series, leaving exactly the product in \eqref{eq:wick}.
	Adding over all pairings proves the formula.
\end{proof}

For example, the fourth moment is
\[
 \mathbb E(\xi_a\xi_b\xi_c\xi_d)=R_{a,b}R_{c,d}+R_{a,c}R_{b,d}+R_{a,d}R_{b,c}.
\]
It follows from \autoref{lem:app-wick} that every fixed Gaussian moment is a polynomial in the
covariance entries, so convergence of covariance matrices implies convergence of all these
moments.

\subsection{Independent rows and uniform moment estimates}\label{app:rows}

We now prove the estimates used in \autoref{thm:clt}.
The elements within a row may fail to commute and the rows may have different distributions.
Only their independence, centering, and a common norm bound are needed for the counting argument.
We first isolate that argument so that it can be applied both to fluctuation elements and to
centered covariance elements.

For a positive integer $d$, denote by $b_d$ the number of partitions of $\{1,\ldots,d\}$.
Its precise value will not matter; in particular, it depends only on the length of the word.
For a subset $D$ of this set, a product indexed by $D$ below is always taken in increasing order.

\begin{lma}\label{lem:app-counting}
Let $N,d\ge1$ be integers, let $C$ be a unital C*-algebra with a tracial state $\tau$, and let
$C_1,\ldots,C_N$ be independent pairwise commuting unital C*-subalgebras of $C$.
Choose constants $K_1,\ldots,K_d\ge0$.
For each $r=1,\ldots,N$ and $a=1,\ldots,d$, let $y_{r,a}\in C_r$ satisfy $\tau(y_{r,a})=0$ and
$\|y_{r,a}\|\le K_a$.
Then
\begin{equation}\label{eq:app-counting}
 \Bigl|\sum_{r_1,\ldots,r_d=1}^N
 \tau(y_{r_1,1}\cdots y_{r_d,d})\Bigr|
 \le b_d N^{\lfloor d/2\rfloor}\prod_{a=1}^d K_a.
\end{equation}
More precisely, the total contribution of tuples inducing a fixed partition with $s$ blocks has
absolute value at most $N^s\prod_{a=1}^d K_a$, and it is zero if that partition has a singleton
block.
\end{lma}
	\begin{proof}
	A tuple $(r_1,\ldots,r_d)$ determines a partition of $\{1,\ldots,d\}$ by placing positions
	$a$ and $b$ in the same block exactly when $r_a=r_b$.
	Fix such a partition $\mathcal P$ and an assignment of distinct row indices to its blocks.
	For $D\in\mathcal P$, denote its assigned row by $r(D)$.
	Commutation between different rows allows us to group the factors by blocks, and independence then gives
	\begin{equation}\label{eq:app-factorization}
 \tau(y_{r_1,1}\cdots y_{r_d,d})
 =\prod_{D\in\mathcal P}\tau\Bigl(\prod_{a\in D}y_{r(D),a}\Bigr).
	\end{equation}
	If $D=\{a\}$ is a singleton, its factor on the right is zero by the assumption $\tau(y_{r(D),a})=0$.
	Otherwise, the absolute value of the right hand side is at most $\prod_{a=1}^d K_a$ by
	contractivity of the state and submultiplicativity of the norm.
	A partition with $s$ blocks has at most $N^s$ assignments of distinct row indices, proving
	the last assertion in the statement.
	Every partition without singleton blocks has at most $\lfloor d/2\rfloor$ blocks.
	There are at most $b_d$ such partitions, so adding their bounds proves \eqref{eq:app-counting}.
\end{proof}

We recall explicitly the normalization to which this will be applied.
Let $(C_r)_{r\in\N}$ be independent pairwise commuting unital subalgebras of $(C,\tau)$, and
let $J$ be a nonempty finite set.
Suppose that $x_{r,g}\in C_r$ is selfadjoint, with $\tau(x_{r,g})=0$, $\tau(x_{r,g}^2)=1$, and
$\|x_{r,g}\|\le M$ for all $r\in\N$ and $g\in J$.
We must have $M\ge1$.
As in \autoref{not:fluctuations}, set
\[
 z_{N,g}=\frac{1}{\sqrt N}\sum_{r=1}^N x_{r,g}
 \quad\text{and}\quad
 k_{N;g,h}=\frac{1}{N}\sum_{r=1}^N
 \frac{x_{r,g}x_{r,h}+x_{r,h}x_{r,g}}{2}
\]
for $N\in\N$ and $g,h\in J$.
We call these the \emph{fluctuation variables} and the \emph{covariance variables}, respectively.
Set $c_{N;g,h}=k_{N;g,h}-\tau(k_{N;g,h})1$, and call these the \emph{centered covariance variables}.
Set $R_N=(R_{N;g,h})_{g,h\in J}$ with $R_{N;g,h}=\tau(k_{N;g,h})$ for all $g,h\in J$.
Traciality makes $R_N$ real symmetric.
If $t_g\in\R$ for $g\in J$, then
\[
 \sum_{g,h\in J}t_gt_hR_{N;g,h}
 =\frac{1}{N}\sum_{r=1}^N
 \tau\Bigl(\bigl(\sum_{g\in J}t_gx_{r,g}\bigr)^2\Bigr)\ge0
\]
Thus $R_N$ is positive semidefinite with diagonal one, and $|R_{N;g,h}|\le1$ by Cauchy--Schwarz.
The fluctuation variables have normalization $\frac{1}{\sqrt N}$, whereas the covariance
variables have normalization $\frac{1}{N}$.
The latter is responsible for the extra decay in part~\ref{item:centered-moments} below.

\begin{prop}\label{prop:moment-estimates}
Let $C$ be a unital C*-algebra with a tracial state $\tau$, and let $(C_r)_{r\in\N}$ be
independent pairwise commuting unital C*-subalgebras of $C$.
Let $J$ be a nonempty finite set, let $M>0$, and let $x_{r,g}\in C_r$ be selfadjoint elements
satisfying $\tau(x_{r,g})=0$, $\tau(x_{r,g}^2)=1$, and $\|x_{r,g}\|\le M$ for all $r\in\N$ and
$g\in J$.
Define $z_{N,g}$ and $k_{N;g,h}$ by \eqref{eq:zk}, set $c_{N;g,h}=k_{N;g,h}-\tau(k_{N;g,h})1$,
and set $R_N^\tau=(\tau(k_{N;g,h}))_{g,h\in J}$.
Let $\xi_N^\tau=(\xi_{N,g}^\tau)_{g\in J}$ be a centered real Gaussian vector with covariance $R_N^\tau$.
Then the following estimates hold.
\begin{enumerate}[(a)]
\item\label{item:gaussian-moments} For every integer $d\ge1$, there is a constant $K_{d,M}>0$, depending only on $d$ and $M$, such that
\begin{equation}\label{eq:moment-error}
 \bigl|\tau(z_{N,g_1}\cdots z_{N,g_d})
       -\mathbb E(\xi_{N,g_1}^\tau\cdots\xi_{N,g_d}^\tau)\bigr|
 \le\frac{K_{d,M}}{\sqrt N}
\end{equation}\nopagebreak
for all $N\in\N$ and $g_1,\ldots,g_d\in J$.
\item\label{item:centered-moments} For all integers $d\ge0$ and $r\ge1$, there is a constant $K_{d,r,M}>0$, depending only on $d$, $r$, and $M$, such that
\begin{equation}\label{eq:centered-error}
 |\tau(w)|\le\frac{K_{d,r,M}}{N^{r/2}}
\end{equation}\nopagebreak
whenever $w$ is a product containing $d$ factors from $(z_{N,g})_{g\in J}$ and $r$ factors from
$(c_{N;g,h})_{g,h\in J}$, in any prescribed order and allowing repetitions, for a fixed
$N\in\N$.
\end{enumerate}
\end{prop}
	\begin{proof}
	We first prove part~\ref{item:gaussian-moments}.
	Fix $d\ge1$ and labels $g_1,\ldots,g_d\in J$.
	Expanding the fluctuation variables gives
	\begin{equation}\label{eq:app-expand}
 \tau(z_{N,g_1}\cdots z_{N,g_d})
 =\frac{1}{N^{d/2}}\sum_{r_1,\ldots,r_d=1}^N
 \tau(x_{r_1,g_1}\cdots x_{r_d,g_d}).
	\end{equation}
	Suppose first that $d$ is odd.
	A partition without singleton blocks has at most $\frac{d-1}{2}$ blocks.
	Applying the counting estimate and dividing by $N^{d/2}$ therefore gives
	\[
 \bigl|\tau(z_{N,g_1}\cdots z_{N,g_d})\bigr|
 \stackrel{\eqref{eq:app-counting}}{\le}\frac{b_dM^dN^{(d-1)/2}}{N^{d/2}}
 =\frac{b_dM^d}{\sqrt N}.
	\]
	The corresponding Gaussian moment is zero by \autoref{lem:app-wick}, proving the assertion for odd $d$.

	Suppose that $d=2m$.
	By \autoref{lem:app-counting}, singleton blocks give zero, and a partition without
	singleton blocks which is not a pairing has at most $m-1$ blocks.
	The total contribution of all these partitions to \eqref{eq:app-expand} therefore has
	absolute value at most $b_dM^d/N$.
	It remains to compare the pairing terms with the Gaussian moment formula.
	Fix a pairing $\mathcal P$, enumerate its pairs as $\{a_k,b_k\}$ with $a_k<b_k$ for
	$k=1,\ldots,m$, and denote its contribution by $A_{N,\mathcal P}$.
	By \eqref{eq:app-factorization}, we have
	\[
 A_{N,\mathcal P}
 =\frac{1}{N^m}
 \sum_{\substack{1\le r_1,\ldots,r_m\le N\\r_1,\ldots,r_m\text{ distinct}}}
 \prod_{k=1}^m\tau(x_{r_k,g_{a_k}}x_{r_k,g_{b_k}}).
	\]
	Let $B_{N,\mathcal P}$ be the same sum without the distinctness condition.
	For $m\ge2$, a tuple excluded by that condition has $r_a=r_b$ for at least one pair of distinct positions $a,b$.
	There are $\binom m2$ choices of these positions and at most $N^{m-1}$ tuples for each choice.
	Each product has absolute value at most $M^{2m}$, so
	\begin{equation}\label{eq:app-collision}
 |A_{N,\mathcal P}-B_{N,\mathcal P}|
 \le\binom m2\frac{M^d}{N}.
	\end{equation}
	The same bound holds for $m=1$, since both sums are then equal.
	The unrestricted sum factors, and traciality identifies its factors with the covariance entries:
	\[
 B_{N,\mathcal P}
 =\prod_{k=1}^m\Bigl(\frac{1}{N}\sum_{r=1}^N
      \tau(x_{r,g_{a_k}}x_{r,g_{b_k}})\Bigr)
 =\prod_{k=1}^m R_{N;g_{a_k},g_{b_k}}.
	\]
	Adding over pairings and using \eqref{eq:wick} now gives
	\[
 \bigl|\tau(z_{N,g_1}\cdots z_{N,g_d})
 -\mathbb E(\xi_{N,g_1}^\tau\cdots\xi_{N,g_d}^\tau)\bigr|
 \le\Bigl(b_d+|\mathcal P_2(d)|\binom m2\Bigr)\frac{M^d}{N}.
	\]
	Since $N\ge1$, this proves part~\ref{item:gaussian-moments}, with
	\[
 K_{d,M}=\Bigl(b_d+|\mathcal P_2(d)|\binom{\lfloor d/2\rfloor}{2}\Bigr)M^d.
	\]
	The count involves only the positions of the word, so repeated labels and the size of $J$ do not affect the bound.

	For part~\ref{item:centered-moments}, the covariance summands must first be centered separately in each row.
	For $s\in\N$ and $g,h\in J$, set
	\[
 v_{s;g,h}=\frac{x_{s,g}x_{s,h}+x_{s,h}x_{s,g}}{2}
 \quad\text{and}\quad
 y_{s;g,h}=v_{s;g,h}-\tau(v_{s;g,h})1.
	\]
	Then $y_{s;g,h}\in C_s$, $\tau(y_{s;g,h})=0$, and $\|y_{s;g,h}\|\le2M^2$.
	Moreover, there is an exact identity
	\begin{equation}\label{eq:app-row-centering}
 k_{N;g,h}-\tau(k_{N;g,h})1
 =\frac{1}{N}\sum_{s=1}^N y_{s;g,h}.
	\end{equation}
	Fix an ordered word $w$ with $d$ fluctuation factors and $r$ centered covariance factors.
	Expanding its sums using \eqref{eq:app-row-centering} produces $d+r$ centered row factors,
	with normalization $\frac{1}{N^{d/2+r}}$.
	The factors have bounds $M$ or $2M^2$ according to their type, regardless of their positions in the word.
	Thus \autoref{lem:app-counting} gives
	\[
 |\tau(w)|
 \le b_{d+r}2^rM^{d+2r}N^{\lfloor(d+r)/2\rfloor-d/2-r}
 \le\frac{b_{d+r}2^rM^{d+2r}}{N^{r/2}}.
	\]
	This proves \eqref{eq:centered-error} with $K_{d,r,M}=b_{d+r}2^rM^{d+2r}$.
\end{proof}

The constants in the proof above are independent of the trace, the labels, and the covariance matrix.
In particular, the trace can vary with $N$, and no lower bound on the eigenvalues of the
covariance matrix is required.
These are the two uniformity properties needed when the central amenability weights are
evaluated at different extreme traces.

\begin{rem}\label{rem:app-replacement}
We explain how the two estimates combine.
Fix a word in fluctuation variables and covariance variables, and replace each covariance
factor $k_{N;g,h}$ by the sum of $\tau(k_{N;g,h})1$ and its centered part.
Expand the word without changing the order of the remaining factors.
The scalar factors have absolute value at most one, and every term containing $r\ge1$ centered
factors tends to zero by \eqref{eq:centered-error}.
The only remaining term has all covariance factors replaced by scalars, and its fluctuation
moment is approximated by \eqref{eq:moment-error}.
Hence, if the covariance entries converge along a sequence, all ordered joint moments converge
to those of a Gaussian vector and its scalar covariance entries.
The detailed application to extreme traces of the row algebra is given in \autoref{thm:clt}.

For example, for a fixed label $g$, the identity $k_{N;g,g}=\frac{1}{N}\sum_{s=1}^Nx_{s,g}^2$ gives
\[
 \tau(z_{N,g}^2)=1
 \quad\text{and}\quad
 \tau((k_{N;g,g}-1)^2)\le\frac{K_{0,2,M}}{N}.
\]
The first element retains a nontrivial fluctuation, whereas the second becomes scalar in tracial $2$-norm.
Part~\ref{item:centered-moments} is stronger than this last observation: it allows the centered
factor to occur between fluctuation factors whose norms may grow with $N$.
\end{rem}

\subsection{From moments to continuous functions}\label{app:moments}

The row estimates concern polynomials, whereas \autoref{sec:gaussian} uses exponentials of
selfadjoint elements whose norms need not be uniformly bounded.
The method of moments gives the passage between these two kinds of functions.
We recall its statements with precise references, rather than reproducing the classical measure-theoretic proofs.
A sequence $(\mu_n)_{n\in\N}$ of probability measures on $\R$ converges \emph{weakly} to $\mu$
if $\int_\R f\,d\mu_n$ converges to $\int_\R f\,d\mu$ for every bounded continuous function
$f\colon\R\to\C$.
A probability measure is \emph{determined by its moments} if it has finite moments of all
orders and no other probability measure has the same moments.

\begin{lma}\label{lem:app-determinate}
Let $\mu$ be a probability measure on $\R$.
If $\int_\R e^{c|x|}\,d\mu(x)<\infty$ for some $c>0$, then all moments of $\mu$ are finite and
$\mu$ is determined by its moments.
\end{lma}
	\begin{proof}
	For $0<r<c$, the power series of $e^{r|x|}$ and monotone convergence give
	\[
 \sum_{k=0}^\infty\frac{r^k}{k!}\int_\R |x|^k\,d\mu(x)
 =\int_\R e^{r|x|}\,d\mu(x)<\infty.
	\]
	Thus the moment power series has positive radius of convergence, so \cite[Theorem 30.1,
	p.~388]{Billingsley} applies.
\end{proof}

\begin{prop}[{\cite[Theorem 30.2, p.~390]{Billingsley}}]\label{prop:app-moment-method}
Let $\mu$ and $(\mu_n)_{n\in\N}$ be probability measures on $\R$ with finite moments of all orders.
Suppose that $\mu$ is determined by its moments and that
\[
 \lim_{n\to\infty}\int_\R x^r\,d\mu_n(x)=\int_\R x^r\,d\mu(x)
\]
for every integer $r\ge0$.
Then $\mu_n$ converges weakly to $\mu$ as $n\to\infty$.
\end{prop}

The assumption that $\mu$ is determined by its moments identifies the possible limit measure, while convergence of the
second moments prevents mass from escaping to infinity.
This is why the statement does not require the supports of the approximating measures to lie in
a common bounded interval.
Recall that a real random variable $X$ is said to have a \emph{finite exponential moment of its
absolute value} if $\mathbb E\bigl(e^{c|X|}\bigr)<\infty$ for some $c>0$.
The following is the C*-algebraic form of the method of moments that we will use.

\begin{cor}\label{cor:app-spectral-convergence}
For $n\in\N$, let $a_n$ be a selfadjoint element of a unital C*-algebra $A_n$, and let $\rho_n$
be a state on $A_n$.
Let $X$ be a real random variable with a finite exponential moment of its absolute value.
If $\lim_{n\to\infty}\rho_n(a_n^r)=\mathbb E(X^r)$ for every integer $r\ge0$, then
\[
 \lim_{n\to\infty}\rho_n(f(a_n))=\mathbb E(f(X))
\]
for every bounded continuous function $f\colon\R\to\C$.
\end{cor}
	\begin{proof}
	Apply \autoref{prop:app-moment-method} to the spectral distributions from
	\autoref{rem:app-spectral} and the law of $X$, which is determined by its moments by
	\autoref{lem:app-determinate}.
\end{proof}

In particular, \autoref{cor:app-spectral-convergence} applies to $f(x)=e^{itx}$ for each fixed $t\in\R$.
Suppose that a selfadjoint polynomial in several selfadjoint elements has a Gaussian quadratic
expression as its limiting model.
Every power of that polynomial is again a polynomial in the original elements, so convergence
of their ordered joint moments gives convergence of all moments of the selfadjoint polynomial.
Once exponential integrability of the limiting quadratic expression is known, the corollary
therefore gives convergence of the traces of its exponentials.
This is the passage used in \autoref{prop:unitaries}; we check the exponential integrability in
the next proposition.

\subsection{Quadratic forms and exponentials}\label{app:quadratic}

Let $V$ be a finite dimensional real Euclidean space, let $\eta$ be a standard Gaussian vector
on $V$, and let $T$ be a selfadjoint operator on $V$.
Set $Q_T=\langle\eta,T\eta\rangle-\Tr(T)$, and call this the \emph{centered Gaussian quadratic
form} associated with $T$.
The covariance of the coordinates of $\eta$ is the identity, so
$\mathbb E\bigl(\langle\eta,T\eta\rangle\bigr)=\Tr(T)$ and $\mathbb E(Q_T)=0$.
When comparing two such quadratic forms, we always use the same vector $\eta$.
We use the unnormalized trace $\Tr$, the Hilbert--Schmidt norm
$\HS T=\Tr(T^*T)^{\tfrac{1}{2}}$, and the trace norm $\|T\|_1=\Tr(|T|)$.

The Fourier transform of the law of $Q_T$ can be estimated using only $\HS T$.
This will give unitaries with small moments and will allow both the covariance and the weights
to vary with the trace.

\Needspace{8\baselineskip}
\begin{prop}\label{prop:gaussianinput}
The following statements hold.
\begin{enumerate}[(a)]
\item\label{item:sqrt} For positive operators $S,T$ on the same finite dimensional real Euclidean space, one has
\[
 \Bigl\|S^{\frac{1}{2}}-T^{\frac{1}{2}}\Bigr\|_{\mathrm{HS}}^2\le\|S-T\|_1.
\]
\item\label{item:quadratic} Let $\eta$ be a standard Gaussian vector on a finite dimensional real Euclidean space, and let $T$ be a selfadjoint operator on that space. Set $Q_T=\langle \eta,T\eta\rangle-\Tr(T)$. Then
\begin{equation}\label{eq:gaussianest}
 \|Q_T\|_{L^2}=\sqrt2\HS T
 \quad\text{and}\quad
 |\mathbb{E}\bigl(e^{isQ_T}\bigr)|\le\Bigl(1+4s^2\HS T^2\Bigr)^{-\frac{1}{4}}
\end{equation}\nopagebreak
for all $s\in\R$.
The law of $Q_T$ has a finite exponential moment of its absolute value and is determined by its moments.
Since $Q_T-Q_U=Q_{T-U}$ when the same $\eta$ is used, the first identity also gives
\begin{equation}\label{eq:gaussiandifference}
 \|Q_T-Q_U\|_{L^2}=\sqrt2\HS{T-U}
\end{equation}\nopagebreak
for all selfadjoint $T,U$ on that space.
\item\label{item:spectral} For each $j\in\N$, let $A_j$ be a unital C*-algebra, let $a_j\in A_j$ be selfadjoint, and let $\rho_j$ be a state on $A_j$.
Let $Q_T$ be a Gaussian quadratic form as in part~\ref{item:quadratic}. If
\[
 \lim_{j\to\infty}\rho_j(a_j^r)=\mathbb{E}(Q_T^r)
\]\nopagebreak
for every integer $r\ge0$, then $\lim_{j\to\infty}\rho_j(f(a_j))=\mathbb E\bigl(f(Q_T)\bigr)$
for every bounded continuous function $f\colon\R\to\C$.
(No uniform bound on $\|a_j\|$ is required.)
\end{enumerate}
\end{prop}
	\begin{proof}
	Part~\ref{item:sqrt} is the finite dimensional Powers--St\o rmer inequality \cite[Lemma 4.1]{PS}.
	Part~\ref{item:quadratic} is immediate in dimension zero, so assume that the Euclidean
	space has dimension $d\ge1$ and let $\lambda_1,\ldots,\lambda_d$ be the eigenvalues of $T$,
	counted with multiplicity.
	In an orthonormal eigenbasis, the coordinates $U_1,\ldots,U_d$ of $\eta$ are independent
	standard real Gaussian variables, and
	\[
 Q_T=\sum_{j=1}^d\lambda_j(U_j^2-1).
	\]
	Each summand has mean zero, and $\mathbb E((U_j^2-1)^2)=2$.
	Independence therefore gives $\mathbb E(Q_T^2)=2\sum_{j=1}^d\lambda_j^2=2\HS T^2$, proving
	the first identity in \eqref{eq:gaussianest}.
	For a standard real Gaussian variable $U$, direct integration gives
	\[
 \mathbb E(e^{zU^2})=\frac{1}{\sqrt{2\pi}}\int_\R e^{-(\frac{1}{2}-z)t^2}\,dt=\frac{1}{\sqrt{1-2z}}
	\]
	for $z\in\C$ with $\operatorname{Re}(z)<\frac12$, taking the square root with positive real part.
	For real $z<\frac12$ this follows by a change of variables in the Gaussian integral, and
	the complex identity follows by analytic continuation on this half plane.
	Applying this identity with $z=is\lambda_j$, and using independence, and the inequality
	$\prod_{j=1}^d(1+a_j)\ge1+\sum_{j=1}^da_j$ for nonnegative $a_1,\ldots,a_d$, we also obtain
	\[
 \bigl|\mathbb E\bigl(e^{isQ_T}\bigr)\bigr|
 =\prod_{j=1}^d\bigl(1+4s^2\lambda_j^2\bigr)^{-\frac{1}{4}}
 \le\Bigl(1+4s^2\sum_{j=1}^d\lambda_j^2\Bigr)^{-\frac{1}{4}}
	\]\nopagebreak
	for every $s\in\R$.
	This is the second estimate in \eqref{eq:gaussianest}.

	Choose $c>0$ such that $2c\max_{1\le j\le d}|\lambda_j|<1$.
	Since $|Q_T|\le\sum_{j=1}^d|\lambda_j|(U_j^2+1)$, another application of the Gaussian integral gives
	\[
 \mathbb E\bigl(e^{c|Q_T|}\bigr)
 \le e^{c\sum_{j=1}^d|\lambda_j|}
       \prod_{j=1}^d\frac{1}{\sqrt{1-2c|\lambda_j|}}<\infty.
	\]
	By \autoref{lem:app-determinate}, the law of $Q_T$ is determined by its moments.
	The identity $Q_T-Q_U=Q_{T-U}$ proves \eqref{eq:gaussiandifference}.

	Part~\ref{item:spectral} follows from \autoref{cor:app-spectral-convergence}, using the
	exponential moment established in part~\ref{item:quadratic}.
\end{proof}

The bounds in \autoref{prop:gaussianinput} explain the two requirements on the quadratic forms
in the main argument.
A lower bound for $\HS T$ makes the nonzero moments of $e^{itQ_T}$ small when $t$ is large,
while a small Hilbert--Schmidt difference between two operators makes their exponentials close
in $L^2$ for a fixed $t$.

\begin{rem}\label{rem:app-square-root}
Let $P$ be a positive operator on a finite dimensional real Euclidean space with $\Tr(P)=1$.
Although there is no dimension independent lower bound for $\HS P$, we always have
$\Bigl\|P^{\tfrac{1}{2}}\Bigr\|_{\mathrm{HS}}=1$, independently of the dimension of the space.
Therefore, \eqref{eq:gaussianest} gives
\[
 \bigl|\mathbb E\bigl(e^{itQ_{P^{\frac{1}{2}}}}\bigr)\bigr|\le(1+4t^2)^{-\frac{1}{4}}
\]
for all $t\in\R$, independently of the dimension and the rank of $P$.
If $P'$ is another positive operator on the same space, then
\[
 \|Q_{P^{\frac{1}{2}}}-Q_{(P')^{\frac{1}{2}}}\|_{L^2}^2
 \stackrel{\eqref{eq:gaussiandifference}}{=}2\Bigl\|P^{\frac{1}{2}}-(P')^{\frac{1}{2}}\Bigr\|_{\mathrm{HS}}^2
 \le2\|P-P'\|_1.
\]
Thus the square root supplies a nonzero variance without losing continuity under the trace norm
perturbations furnished by amenability.
In \autoref{sec:gaussian}, the operators $P$ are positive weighted sums of rank one operators
associated to unit Gram vectors.
Their trace is one because the weights add up to one.
The estimate therefore holds even when the Gram vectors are linearly dependent.

\end{rem}

We finish by recording the perturbation estimate for exponentials on the C*-algebra side.
It avoids any need to compare joint spectral measures of two noncommuting elements.

\begin{lma}\label{lem:app-duhamel}
Let $C$ be a unital C*-algebra, let $\tau\in T(C)$, and let $a,b\in C$ be selfadjoint.
Then
\begin{equation}\label{eq:duhamel}
 \|e^{ita}-e^{itb}\|_{2,\tau}\le |t|\|a-b\|_{2,\tau}
\end{equation}
for all $t\in\R$.
\end{lma}
	\begin{proof}
	For $t\ge0$, define $f\colon[0,t]\to C$ by $f(u)=e^{i(t-u)a}e^{iub}$ for $u\in[0,t]$.
	This function is norm differentiable, with $f'(u)=i e^{i(t-u)a}(b-a)e^{iub}$, so integration yields
	\[
 e^{ita}-e^{itb}=i\int_0^t e^{i(t-u)a}(a-b)e^{iub}\,du.
	\]
	Left and right multiplication by unitaries preserve the tracial $2$-seminorm.
	Taking that seminorm and using the triangle inequality for the integral proves \eqref{eq:duhamel} for $t\ge0$.
	For $t<0$, apply the result to $-a$ and $-b$ with the positive parameter $-t$.
\end{proof}

Suppose now that a sequence of selfadjoint elements has moments converging to those of $Q_T$.
By \autoref{cor:app-spectral-convergence}, the traces of their exponentials converge to
$\mathbb E\bigl(e^{itQ_T}\bigr)$, and \eqref{eq:gaussianest} bounds the absolute value of this
limit.
For two sequences whose joint moments converge to two quadratic forms on the same Gaussian
space, the trace of the square of their difference converges to the squared $L^2$-norm of the
difference of the forms.
Equations \eqref{eq:gaussiandifference} and \eqref{eq:duhamel} then bound the difference of their exponentials.
These statements apply along sequences of varying traces because the estimates in
\autoref{prop:moment-estimates} are uniform in the trace.
In the proof of \autoref{prop:unitaries}, compactness of the possible covariance matrices and
scalar weights allows us to extract such a sequence from any proposed failure of a uniform
estimate.
This is how the scalar Gaussian calculations yield the uniform tracial estimates required for
local equivariant property Gamma.

\end{document}